\documentclass[11pt]{article}
\usepackage{graphicx, overpic, fancyhdr}
\usepackage{url}
\usepackage{color}
\usepackage{latexsym,array,delarray,amsmath,amsfonts,amssymb,amsthm}
\usepackage{natbib}
\usepackage{bm}
\bmdefine{\Bt}{t}
\bmdefine{\BX}{X}
\bmdefine{\BY}{Y}
\bmdefine{\BZ}{Z}
\bmdefine{\BB}{B}
\bmdefine{\BM}{M}
\bmdefine{\BD}{D}
\bmdefine{\Bi}{i}
\bmdefine{\Bj}{j}
\bmdefine{\Bx}{x}
\bmdefine{\By}{y}
\bmdefine{\Bz}{z}
\bmdefine{\Bw}{w}
\bmdefine{\Ba}{a}
\bmdefine{\Bb}{b}
\bmdefine{\Bc}{c}
\bmdefine{\Be}{e}
\bmdefine{\Bh}{h}
\bmdefine{\Bg}{g}
\bmdefine{\Bu}{u}
\bmdefine{\Bv}{v}

\def\Qsat{Q_{\rm sat}}

\theoremstyle{plain}
\newtheorem{thm}{Theorem}[section]
\newtheorem{lemma}[thm]{Lemma}
\newtheorem{prop}[thm]{Proposition}
\newtheorem{cor}[thm]{Corollary}

\newtheorem{defn}[thm]{Definition}
\newtheorem{ex}[thm]{Example}
\newtheorem{pr}[thm]{Problem}
\newtheorem{rmk}[thm]{Remark}

\newcommand{\Q}{{\mathbb Q}}
\newcommand{\Z}{{\mathbb Z}}
\newcommand{\N}{{\mathbb N}}

\newcommand{\C}{{\mathbb C}}
\newcommand{\R}{\mathbb R}

\newcommand{\cone}{{\rm cone}}

\def\comment#1{\textit{[#1]}}
\def\comment#1{}

\begin{document}
\renewcommand{\baselinestretch}{1.2}
%\lhead[\fancyplain{} \leftmark]{}
%\chead[]{}
%\rhead[]{\fancyplain{}\rightmark}
%\cfoot{}
%\headrulewidth=0pt
\markright{
%\hbox{\footnotesize\rm Statistica Sinica
%{\footnotesize\bf ??}(200?), 000-000}\hfill
}
%\pagestyle{myheadings}
%\markboth{\hfill{\footnotesize\rm AKIMICHI TAKEMURA AND RURIKO YOSHIDA
%}\hfill}{\hfill {\footnotesize\rm SATURATION POINTS} \hfill}

\comment{
%\renewcommand{\thefootnote}{}
%$\ $\par
%\fontsize{10.95}{14pt plus.8pt minus .6pt}\selectfont
\vspace{0.8pc}
\centerline{\large\bf Holes in semigroups and the multi-dimensional 
Frobenius problem}
%\vspace{2pt}
%\centerline{\large\bf semigroup and their applications to contingency tables}
\vspace{.4cm}
\centerline{Akimichi Takemura$^*$ and Ruriko Yoshida$^+$}
\vspace{.4cm}
\centerline{\it University of Tokyo$^*$  University of Kentucky$^+$}
\vspace{.55cm}
\fontsize{9}{11.5pt plus.8pt minus .6pt}\selectfont
}

\title{A generalization of the integer linear infeasibility problem}
%\title{Holes in semigroups and the multi-dimensional Frobenius problem}

\author{ Akimichi Takemura  \\ Department of Mathematical Informatics \\ 
   University of Tokyo\\ Bunkyo, Tokyo, Japan \\
   takemura@stat.t.u-tokyo.ac.jp \and Ruriko
  Yoshida \\ Department of Statistics \\ University of Kentucky \\ Lexington,
  KY USA\\ ruriko@ms.uky.edu}

\maketitle

\begin{abstract}
Does a given system of linear equations $A\Bx = \Bb$ have 
a nonnegative integer solution?  
This is a fundamental question in many areas, 
such as operations research, number theory, and statistics.
In terms of optimization, this is called an {\em integer feasibility problem}.
A {\em generalized integer feasibility problem} is to find %all right-hand-side
$\Bb$ such that there does not exist a nonnegative integral solution in the 
system with a given $A$.
One such problem is the well-known {\em Frobenius problem}.  
In this paper we study the generalized integer feasibility problem and 
also the multi-dimensional Frobenius problem. 
%This is a fundamental question in many areas.
%In terms of optimization, this is called an {\em integer feasibility problem}.
%In number theory, we can apply to {\em Frobenius problems}.
%In statistics this problem arises in data security problems for
%contingency table data.% and also is closely related to nonsquarefree
%elements of Markov bases for sampling 
%contingency tables with given marginals.
%In this paper we study a family of systems with the right-hand-side 
%allowed to vary.  This generalization of the integer linear infeasibility 
%problem finds applications in solving {\em Frobenius problems} and
%{\em multi-dimensional integer planar transportation problems}.
To study a family of systems with no nonnegative integer solution, we focus 
on a 
% commutative semigroup generated by the columns of its defining matrix. 
% In this paper we study 
commutative semigroup generated by a finite 
subset of $\Z^d$ and its saturation.  
An element in the difference of the semigroup and its saturation is
called a ``hole''.
We show the necessary and sufficient
conditions for  the finiteness of the set of holes. 
% the given semigroup to have a finite number of elements in 
% the difference between the semigroup and its saturation. 
Also we define fundamental holes and saturation points of a 
commutative semigroup.  Then, we show the simultaneous 
finiteness of the set of holes, 
%difference between the semigroup and its saturation, 
the set of non-saturation points, %of the semigroup, 
and the set of  generators for saturation points. %, which is a set of generator of a monoid.
As examples we consider some three- and four-way contingency tables
  from statistics and apply our results to them.
Then we will discuss the time complexities of our algorithms.
%We apply our results to some three- and four-way contingency tables.
%some statistical data security problems.
\end{abstract}

\vspace{9pt}
\begin{quotation}
\noindent {\it Key words and phrases:}
%Keywords: 
contingency tables, 
data security,
Frobenius problem, 
indispensable move,
Markov basis,
monoid, 
Hilbert basis, 
linear integer feasibility problem, 
saturation, 
semigroup 
\par
\end{quotation}\par

\par
\fontsize{10.95}{14pt plus.8pt minus .6pt}\selectfont
% \setcounter{chapter}{1}
% \setcounter{equation}{0} %-1
%\noindent {\bf 1. Introduction}
\section{Introduction}

Consider the following system of linear equations and inequalities:
\begin{equation}\label{feas}
A\Bx = \Bb, \, \, \, \Bx \geq 0, 
\end{equation}
where $A \in \Z^{d\times n}$ and $\Bb \in \Z^d$.
Suppose the solution set $\{x \in \R^n: A\Bx = \Bb, \, \Bx \geq 0\} \not
= \emptyset$.  
The {\em linear integer feasibility problem} is to ask whether the system in 
\eqref{feas}
has an integral solution or not.  
A {\em generalized integer feasibility problem} is to find all %right-hand-side
$\Bb$ such that there does not exist a nonnegative integral solution in the 
system with a given $A$.
Note that there exists an integral solution
for the system in \eqref{feas} if and only if $\Bb$ is in the {\em semigroup}
generated by the column vectors of $A$.  From this, 
we can write this problem as follows.
\begin{pr}\label{problem1}
Let $\Ba_1,\ldots,\Ba_n \in \Z^d$ be columns of $A$ and
\begin{equation}
\label{eq:semigroupQ}
Q = Q(A)= \left\{\Ba_1x_1 + \cdots + \Ba_n x_n : x_1, \dots, x_n \in   \Z_+ \right\}
\end{equation}
be the set of all nonnegative integer combinations of $\Ba_1,\ldots,\Ba_n$
or in other words the {\em semigroup} % $Q \subset \Z^d$ 
generated by $\Ba_1,\ldots,\Ba_n$.  Compute a finite representation 
of all vectors of $Q$.
\end{pr}
%Thus solving Problem \ref{problem1}
%means to compute a finite representation of all elements in the 
%difference between the semigroup generated by the columns of $A$ and its 
%{\em saturation}. 
\cite{Barvinok2002} introduced an algorithm to encode all 
vectors in the semigroup $Q$ 
into a generating function as a {\em short rational generating function}
in polynomial time when $d$ and $n$ are fixed. Therefore, using this algorithm
one can compute a finite representation 
of all vectors not in $Q$ in polynomial time if we fix $d$ and $n$.
However, their algorithm is 
yet technically difficult to implement so that we do not know whether it is 
practical or not.  
Modifying Problem \ref{problem1}, 
in this paper, we would like to solve
the following problem.  
\begin{pr}\label{problem2}
% Let $\Ba_1,\ldots,\Ba_n \in \Z^d$ be columns of $A$ and
% \[
% Q = \left\{\Ba_1x_1 + \cdots + \Ba_n x_n: x_n \in \Z_+^d\right\}
% \]
% be the set of all nonnegative integer combinations of $\Ba_1,\ldots,\Ba_n$
% or in other words the {\em semigroup} $Q \subset \Z_+^d$ 
% generated by $\Ba_1,\ldots,\Ba_n$.  
Let $Q$ be defined as in Problem \ref{problem1}.
Decide whether there is a finite number 
of integral vectors not in $Q$ but in its saturation.
\end{pr}
In other words, for 
fixed $A$, decide whether there is a finite number of integral vectors 
$\Bb \in \Z^d$ 
such that the system in \eqref{feas} has a nonnegative rational solution 
but not a nonnegative integral solution.

Intensive research has been carried out on integer feasibility problems.
%For a {\em fixed vector $\Bb$},
In 1972, \cite{karp} showed that solving the integer linear feasibility
problem is NP hard.
In the 1980's, H.W. Lenstra, Jr. developed an algorithm to {\em detect}
integer solutions in the system \eqref{feas} using the LLL-algorithm
[\cite{grolosch,lenstra}].  Lenstra also showed that
integer programming problems with a fixed number of variables can be
solved in time polynomial in the input size. The algorithm was actually
developed in order to prove that the integer feasibility problem can
be solved in polynomial time if the dimension is fixed. A later algorithm of
similar structure, by \cite{lovaszscarf}, was
implemented by \cite{cooketal}. In addition, Aardal and
collaborators [\cite{aardaletal3,aardaletal2,aardaletal1}] have
used the LLL-procedure to rewrite a system of linear equations 
into an equivalent system that was easier
to solve with the branch-and-bound method 
for testing integer feasibility. In the 1990's, based on work by the
geometers Brion, Khovanski, Lawrence, and Pukhlikov, Barvinok
discovered an algorithm to {\it count} integer points in rational
polytopes, and this algorithm also runs in polynomial time if we fix
the dimension [\cite{bar,BarviPom}]. The idea of the algorithm is to
encode all the integer solutions for the system in \eqref{feas} into a
rational generating function. 

In recent years, the generalized integer linear feasibility problem has found 
applications 
in many research areas, such as number theory and statistics. 
%In number theory, the {\em Frobenius problem} is to find the biggest positive
% integer $b$ such that there does not exist integral solution in \eqref{feas}
%with $d = 1$ [\cite{aardaletal3}].
One such problem is 
the well-known {\em Frobenius problem}, that is, for $d=1$ 
and relatively prime positive integers $a_1, \dots, a_n$,
it is to find the biggest positive
integer $b$ such that there does not exist an 
integral solution in \eqref{feas}
[\cite{aardaletal3}].  Equivalently, it is to find the
smallest positive integer $b'$ such that there exists an integral 
solution with $b = b' + \bar{b}$ for any $\bar{b} \in \Z_+$ in \eqref{feas}. 
\comment{More precisely, we have the following:
Let $\Ba_1,\ldots,\Ba_n$ be positive coprime integers and
\[
Q = \left\{\Ba_1x_1 + \cdots + \Ba_n x_n: x_n \in \Z_+\right\}
\]
be the set of all nonnegative integer combinations of $\Ba_1,\ldots,\Ba_n$
or in other words the {\em semigroup} $Q \subset \Z_+$ of nonnegative integers
generated by $\Ba_1,\ldots,\Ba_n$.  The Frobenius problem is to 
find the smallest positive integer $b' \in Q$ such that 
$\Bb = b' + \bar{b}$ for any $\bar{b} \in \Z_+$ in $Q$ and 
how many positive integers are not in $Q$.
}
Since Georg Frobenius focused  on this problem, 
it attracted substantial attention over more than a hundred years
(see [\cite{frob}] for a nice survey).
We can generalize the Frobenius problem to the multi-dimensional
case.  Let $\Ba_1,\ldots,\Ba_n \in \Z^d$ such that 
the lattice $L$ generated by them is $\Z^d$. 
Let $K=\cone(\Ba_1,\allowbreak\ldots,\Ba_n)$ be the cone generated by 
$\Ba_1,\ldots,\Ba_n$ 
and let
% \[
% Q = \left\{\Ba_1x_1 + \cdots + \Ba_n x_n: x_n \in \Z_+^d\right\}
% \]
$Q$ in (\ref{eq:semigroupQ})
% be the set of all nonnegative integer combinations of $\Ba_1,\ldots,\Ba_n$
% or in other words the {\em semigroup} $Q \subset \Z_+^d$ 
be the semigroup generated by $\Ba_1,\ldots,\Ba_n$. 
Let $S = \{\Bb \in Q: \Bb + (K \cap L) \subset Q\}$.  
In [\cite{Sturmfels2004}], a vector $\Bb \in S$  
is called a {\em saturation point} in $Q$.
% Then, the multi-dimensional Frobenius problem is to
% find all vectors $\Bb \in S$ such that
% (1) $\Bb \in S$ cannot be written as a nonnegative integer combination of elements in 
% $S$, (2) $\Bb \in S$ cannot be written as a nonnegative integer 
% combination of elements in $Q$ or (3) $\Bb \in S$ cannot be written 
% as a nonnegative integer combination of elements in $K\cap L$, and how many 
% vectors are not in $Q$ but in its saturation.
We ask to find ``minimal'' elements of $S$.  In the multi-dimensional
version of the Frobenius problem, 
the notion of minimality can be defined in several ways.
We present three definitions of
minimality and show finiteness results of the set of the 
minimal elements of $S$ for each definition.

% Therefore in this paper, we study the set of all saturation points 
% $\Bb \in Q$ such that $\Bb + (K \cap L) \subset Q$.  We will show that 
% sets of saturation points satisfying (2) and (3) are always finite. 
% Then we will show the simultaneous 
% finiteness of the difference between the semigroup and its saturation, 
% the set of non-saturation points, %of the semigroup, 
% and the set of saturation points satisfying (1).

In statistics, one can find an application in the data security problem of 
{\em multi-way contingency tables} [\cite{dobra-karr-sanil2003}].
%, or {\em tabular data}.
The {\em 3-dimensional integer planar transportation problem} (3-DIPTP)
is an integer feasibility problem which asks whether there exists
a three dimensional contingency table with the given 2-marginals or not. 
(In graph theory, a graph is called 
{\em planar} if it can be drawn in a plane without graph edges crossing.)
For more details on the 3-DIPTP, see [\cite{cox02}].
%Vlach in [Vlach, 1986] 
% 3-DIPTP  
%
\comment{
It has received recent attention in data 
security problem (\cite{Cox02}). %[Cox, 2002].   
In order to publish the data to public, national 
statistical offices (NSOs) subject 
statistical data to a range of verification and ``cleaning'' processes. 
%Much data is derived from probability samples, requiring estimation of 
%population-level statistics. 
Data in a statistical database can come from
multiple sources, at various times, and may have been modified through a 
variety of statistical procedures such as rounding. 
%In 2000, Cox 
\cite{Cox00} demonstrates that any of these factors can produce an 
infeasible table.  Even though marginal totals satisfy obvious necessary 
conditions,  for some instances, there does not exist any feasible solution 
with the marginals.} % (\cite{Cox00}). %[Cox, 2000].
\cite{Vlach1986} provides an excellent summary of attempts on 3-DIPTP.

The linear integer feasibility problem is also closely related to the
theory of Markov bases [\cite{Sturmfels1998}] for sampling contingency
tables with given marginals by Markov chain Monte Carlo methods.
The notion of indispensable moves of Markov bases was defined in 
[\cite{takemura-aoki-2004aism}] and further studied in 
[\cite{ohsugi-hibi-indispensable}].
Recently
\cite{ohsugi-hibi-contingency-tables-2006} gave a simple explicit method
to construct infeasible equations of (\ref{feas}) from
non-squarefree indispensable moves of Markov bases.  %See the
One finds more details in a 
discussion of three-way tables in Section \ref{sec:contingency-table}.

In Section \ref{intro}
we define saturation points and then
we will state our main theorem, Theorem \ref{condition}, which shows
the simultaneous finiteness of the set of {\em holes}, which is the difference 
between the semigroup and its saturation,
the set of {\em non-saturation points} of
the semigroup, and the set of generators for saturation points.  In Section
\ref{necsuf}, we show the necessary and sufficient condition for the
finiteness of the set of holes. % which is the difference 
%between the semigroup and its saturation.  
Section \ref{satpt} shows
% conditions for 
%the simultaneous finiteness of the set of holes, 
%the set of non-saturation points of
%the semigroup, and the set of generators for saturation points.
a proof of Theorem \ref{condition}.
Section \ref{sec:contingency-table} contains 
% some computational simulations on 3-DIPTP.
various computational results for three- and four-way contingency
tables.  Section \ref{disc} will discuss that
(1) solving Problem \ref{problem1},  (2) solving Problem \ref{problem2},
(3) computing the set of {\em holes}, and 
(4) computing the set of {\em fundamental holes} 
are polynomial time in fixed $d$ and $n$.

\section{Notation and the main theorem}\label{intro}

In this section we will remind the reader of some definitions and we will set 
appropriate notation.
We follow the notation in Chapter 7 of [\cite{Sturmfels2004}] and
[\cite{Sturmfels1996}].
%[Miller and Sturmfels, 2004] and 
%[Sturmfels, 1996].  
Let $A \in \Z^{d\times n}$ and let $\Ba_1,\ldots,\Ba_n$ denote the
columns of $A$.  Let $\N=\Z_+=\{0,1,\ldots\}$.  

\begin{defn}
Let $Q$ in (\ref{eq:semigroupQ}) be the semigroup generated by $\Ba_1,\ldots,\Ba_n$, let
$K=\cone(\Ba_1,\allowbreak\ldots,\Ba_n)$ be the cone generated by $\Ba_1,\ldots,\Ba_n$, and
let $L$ be the lattice generated by $\Ba_1,\ldots,\Ba_n$.  
Then the semigroup $\Qsat = K \cap L$ is called the
{\em saturation} of the semigroup $Q$. $Q \subset \Qsat$ and
we call $Q$ {\em saturated} if  $Q = \Qsat$ (also this is called {\em normal}).
$H=\Qsat \setminus Q$ is
the set of holes.  $\Ba \in Q$ is called a saturation point if 
$\Ba + \Qsat \subset Q$. 
\end{defn}
We assume $L=\Z^d$ without loss of generality for our theoretical
developments in Sections 3 and 4.  This is for convenience in working with
the Hilbert basis of $K$.
The following is a list of some notations through this paper:
\begin{eqnarray*}
% Q&=&A {\mathbb N}^n=\{ \Ba_1 x_1  + \cdots +  \Ba_n x_n : 
% \lambda_1, \dots, \lambda_n \in \N \} \\
K&=&A {\mathbb R}_+^n= \{ \Ba_1 x_1 + \cdots +  \Ba_n x_n: 
x_1, \dots, x_n \in \R_+ \}\\
% L&=&A {\mathbb Z}^n =\{ \Ba_1 x_1 + \cdots + \Ba_n x_n: 
% \lambda_1, \cdots, \lambda_n \in \Z \} \\
\Qsat&=&K \cap L = \mbox{saturation of $A$} \supset Q \\
H &=& \Qsat \setminus Q = \mbox{holes in } \Qsat\\
S &=& \{ \Ba \in Q : \Ba + \Qsat \subset Q \}= \mbox{saturation
  points of } Q\\
\bar S &=& Q\setminus S = \mbox{non-saturation
  points of } Q
\end{eqnarray*}
We assume that there
exists $\Bc \in \mathbb{Q}^d$ such that $\Bc \cdot \Ba_i > 0$ for
$i=1,\ldots,n$, where $\cdot$ is the standard inner product.  
% Multiplying $\Bc$ by the least common multiple of the denominators of its
% elements, we can without loss of generality assume that $\Bc \in
% \mathbb{Z}^d$.
Under this assumption $K$ and $Q$ are {\em pointed} and 
$S$ is non-empty by Problem 7.15 of [\cite{Sturmfels2004}]. %[Miller and Sturmfels, 2004].
%Let $\bar S=Q-S$ denote the set of non-saturation points of $Q$.  Then
$\Qsat$ is partitioned as
\[
\Qsat =  H \cup \bar S \cup S  = H \cup Q.
\]
Equivalently  
\begin{equation}
\label{eq:1}
S \subset Q \subset \Qsat 
\end{equation}
and the differences of these two inclusions are $\bar S$ and $H$,
respectively.

If $Q$ is saturated (equivalently
$H=\emptyset$), then $0 \in S$ and $S=Q$, because
$Q=0+Q \subset S+Q \subset S$.
%, and by the monotonicity of $S$, $S=Q$.  
Therefore $S=Q=\Qsat$ in
\eqref{eq:1}. Similarly if $S=Q$, then $0\in S$ and $\Qsat
\subset Q$, implying $Q$ is saturated.  From this consideration it
follows that either $S=Q=\Qsat$ or the two inclusions in \eqref{eq:1}
are simultaneously strict.  
% [careful on the defition of $S$, whether it contains the origin or not]
% Also note that $0\in S$ if and only if $Q$ is saturated.

\comment{
For $\Ba, \Bb \in L$, write
\begin{eqnarray*}
&&\Ba \preceq_Q \Bb \quad     \Leftrightarrow \quad  \Bb - \Ba \in Q \\
&&\Ba \preceq_S \Bb \quad     \Leftrightarrow \quad  \Bb - \Ba \in S
\cup \{0\}\\
&&\Ba \preceq_{Q_{\rm sat}} \Bb \quad     \Leftrightarrow \quad  \Bb -
\Ba \in \Qsat 
\end{eqnarray*}
We also write $\Bb\succeq_Q \Ba$, \ $\Bb\succeq_S \Ba$,
$\Bb\succeq_{\Qsat} \Ba$.  If $\Ba \neq \Bb$, then we use $\prec$
instead of $\preceq$.
\begin{defn}
$\Ba\in S$ is called %a {\em saturation point} if 
\begin{itemize}
\setlength{\itemsep}{0pt}
\item[a)] a $Q$-minimal saturation point if $\Bb \prec_Q \Ba \Rightarrow
\Bb\not\in S$, 
\item[b)] a $S$-minimal saturation point if $\Bb \prec_S \Ba \Rightarrow
\Bb\not\in S$,
\item[c)] a $\Qsat$-minimal saturation point if $\Bb \prec_{\Qsat} \Ba
\Rightarrow
\Bb\not\in S$.
\end{itemize}
\end{defn}
}

\medskip
We now consider three different notions of the minimality of
saturation points, i.e.,  
points of $S$ which are minimal 
with respect to $S$, $Q$, and
$\Qsat$.  We call $\Ba\in S$ an $S$-minimal (a $Q$-minimal,
a $\Qsat$-minimal, resp.) if there exists no other $\Bb\in S$, $\Bb\neq \Ba$,
such that $\Ba- \Bb \in S$ ($Q$, $\Qsat$, resp.).  More formally
$\Ba\in S$ is 
%Denoting $S_0 = S \cup \{0\}$, these definitions can be rewritten as
%Denoting $S_0 = S \setminus \{0\}$, these definitions can be rewritten as
\begin{itemize}
\setlength{\itemsep}{0pt}
\item[a)] an $S$-minimal saturation point if $(\Ba +
  (-(S\cup\{0\})))\cap S=\{ \Ba\}$,
\item[b)] a $Q$-minimal saturation point if $(\Ba + (-Q))\cap S =\{\Ba\}$,
\item[c)] a $Q_{\rm sat}$-minimal saturation point if $(\Ba +
  (-\Qsat))\cap S=\{\Ba\}$.
\end{itemize}

Let $\min(S;S)$ denote the set of $S$-minimal saturation points,  
$\min(S;Q)$ the set of $Q$-minimal saturation points, and 
$\min(S;Q_{\rm sat})$ the set of $Q_{\rm sat}$-minimal saturation points. 
Because of the inclusion \eqref{eq:1}, it follows that
\begin{equation}
\label{eq:2}
\min(S;\Qsat) \subset \min(S;Q)\subset \min(S;S).
\end{equation}

If $\Ba \in H$, then for any $\Bb\in Q$, either $\Ba-\Bb \not\in
\Qsat$ or $\Ba-\Bb \in H$.  This is because if $\Ba-\Bb\in
\Qsat$ and $\Ba-\Bb \not\in H$, then $\Ba-\Bb \in Q$, and hence
$\Ba=\Bb + (\Ba-\Bb) \in Q$, which contradicts $\Ba \in H$.  This relation
can be expressed as
\[
\Qsat \cap (H + (-Q)) = H.
\]
This relation suggests the following definition.

\begin{defn}
We call $\Ba \in \Qsat$, $\Ba\neq 0$, a {\em fundamental hole} if
\[
\Qsat \cap (\Ba + (-Q)) = \{ \Ba \}.
\]
Let $H_0$ be the set of fundamental holes.
\end{defn}

\begin{ex}
Consider the one-dimensional example $A=(3 \ 5 \ 7)$ with $L = \Z$.
$\Qsat=\{0,1,\ldots \}$, $Q=\{0,3,5,6,7,\ldots\}$,
$-Q=\{0,-3,-5,-6,-7,\ldots\}$,   $H=\{1,2,4\}$, 
$S = \{5, 6, 7, \dots\}$ and $\bar{S} = \{0, 3\}$. 
Among the 3 holes, $1$ and $2$
are fundamental.  
For example, $2\in H$ is fundamental because
\[
 \{0,1,\ldots\} \cap \{2,-1,-3,-4,-5,\ldots\} = \{ 2 \}.
\]
On the other hand $4\in H$ is not fundamental because
\[
 \{0,1,\ldots \} \cap \{4,1,-1,-2,-3,\ldots\} = \{4, 1\}.
\]
\end{ex}

If $0\neq \Ba \in Q$, then $\Qsat \cap (\Ba + (-Q)) \supset \{ \Ba, 0 \}$
and $\Ba$ is not a fundamental hole.  This implies that a fundamental
hole is a hole.  For every  non-fundamental hole $\Bx$, there exists
$\By\in H$ such that $0\neq \Bx-\By \in Q$.  If $\By$ is not
fundamental we can repeat this procedure.  Since the procedure has to
stop in finite number of steps, it follows that every non-fundamental
hole $\Bx$ can be written as
\begin{equation}\label{fundholes}
\Bx = \By + \Ba,   \qquad \By \in H_0, \quad \Ba \in Q, \ \Ba\neq 0.
\end{equation}

We also focus on a {\em Hilbert basis} of a cone $K$ and in the next
section we will show
a relation between the set of holes $H$ and the {\em minimal Hilbert basis} of a 
pointed cone $K$. 
\begin{defn}
We call a finite % integral vector 
subset $B \subset K \cap \Z^d$ a {\em Hilbert
basis} of a cone $K$ if any integral point in $K$ can be written as a 
nonnegative integral linear combination of elements in $B$.  If $B$ is 
minimal in terms of inclusion then we call it a {\em minimal Hilbert basis} of
$K$.
\end{defn}
Note that there exists a Hilbert basis for any rational polyhedral cone and 
also if a cone
is pointed then there exists a unique minimal Hilbert basis 
%([Schrijver, 1986] for more details).
[see \cite{Schrijver1986} for more details].
%\vskip 0.3in
\bigskip

Now we will present our main theorem of this paper and then we will present
small examples to demonstrate the theorem.  In the theorem,  $\cone(S)$ 
denotes the set of finite nonnegative real combinations of
elements of $S$ and ``rational polyhedral cone'' is a closed cone defined
by rational linear weak inequalities (inequalities that permit the equality 
case).
One can find a proof of this theorem in Section \ref{satpt}.

\begin{thm}\label{condition}
%If Condition $1$ holds, then $H, \bar S$ and
%$\min(S;S)$ are all finite.  If Condition $1$ does not hold, then $H, \bar
%S$ and $\min(S;S)$ are all infinite.
%$H, \bar S$, and $\min(S;S)$ are simultaneously finite or simultaneously
%infinite.
The following statements are equivalent.  \begin{enumerate}
\item\label{equiv1} $\min(S;S)$ is finite. 
\item\label{equiv2} $\cone(S)$ is a rational polyhedral cone.
\item\label{equiv3} There is some $s\in S$ on every extreme ray of $K$.
\item\label{equiv4} $H$ is finite.
\item\label{equiv5} $\bar S$ is finite.
\end{enumerate}
\end{thm}

\begin{figure}[ht]
\begin{center}
     \includegraphics[scale=.4]{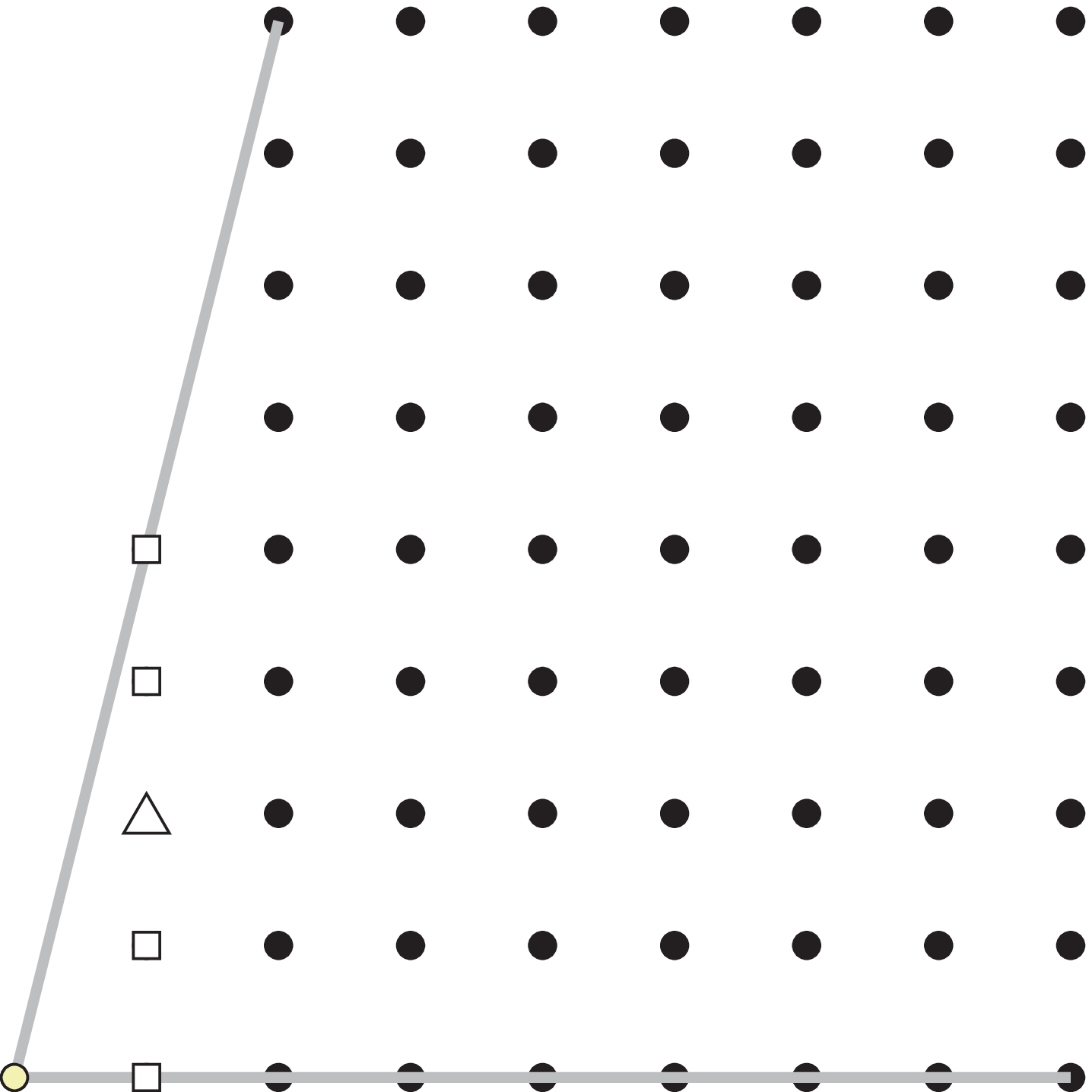}
 \end{center}
\caption{White circles represent non-saturation points, a triangle represents a hole, 
white squares represent $S$-minimal saturation points, and black circles 
represent non $S$-minimal saturation points  
in the semigroup in Example \ref{example1}.}
\end{figure}
\begin{ex}\label{example1}
Let $A$ be an integral matrix such that
\[
A = \left(\begin{array}{cccc}
1 & 1 & 1 & 1\\
0 & 1 & 3 & 4\\
\end{array}\right).
\]
The set of holes $H$ consists of only one element $\{(1, 2)^t\}$.
$\bar S = \{(0, 0)^t\}$. $\min(S;S) = \{
(1, 0)^t, \, (1, 1)^t, \\ (1, 3)^t, \, 
(1, 4)^t\}$.
Thus, $H$, $\bar S$, and $\min(S;S)$ are all finite.
\end{ex}
\begin{figure}[ht]
\begin{center}
     \includegraphics[scale=.4]{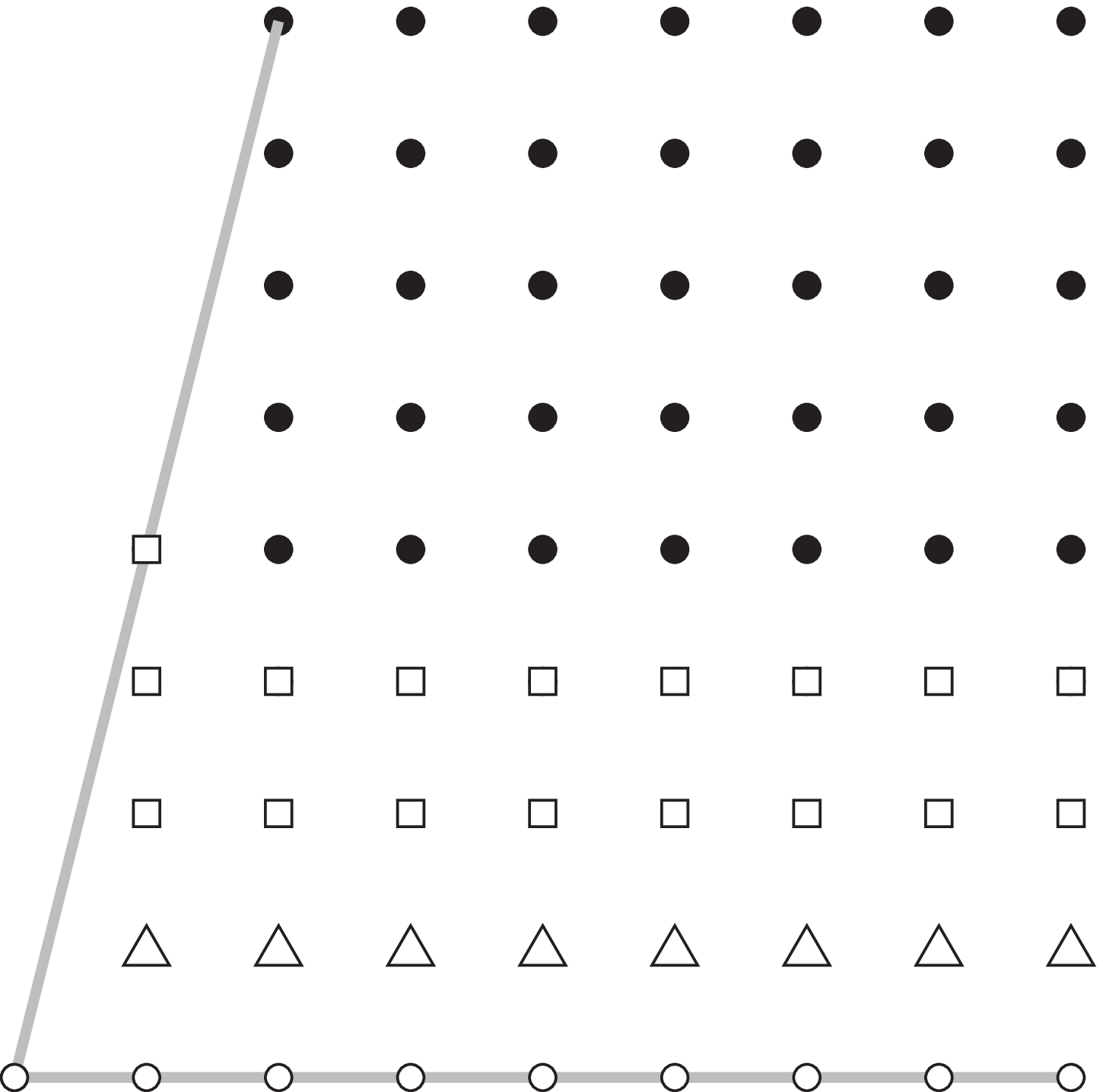}
 \end{center}
\caption{White circles represent non-saturation points, triangles represent holes, 
white squares represent $S$-minimal saturation points, and black circles 
represent non $S$-minimal saturation points  
in the semigroup in Example \ref{example2}.}
\end{figure}
\begin{ex}\label{example2}
Let $A$ be an integral matrix such that
\[
A = \left(\begin{array}{cccc}
1 & 1 & 1 & 1\\
0 & 2 & 3 & 4\\
\end{array}\right).
\]
The set of holes $H$ are the elements $\{(k, 1): k \in \Z , \, 
k \geq 1 \}$.
%$\bar S = \{(i, 0), \, (j, 1): i, \, j \in \Z , \, i, j \geq 1\}$, 
$\bar S = \{(i, 0)^t : i  \in \Z , \, i\geq 0\}$, 
and 
%$\min(S;S)=\{(k, 2): k \in \Z , \, k \geq 1 \}$.
$\min(S;S)=\{(k, j)^t: k \in \Z , \, k \geq 1, \, 
2 \leq j \leq 3 \} \cup \{(1, 4)\}$.
Thus, $H$, $\bar S$, and $\min(S;S)$ are all infinite.
However, $\min(S;Q) = \{(1, 2)^t, \, (1, 3)^t, \, (1, 4)^t\}$ is finite.
\end{ex}

\section{Necessary and sufficient condition of finiteness of a set of
holes}\label{necsuf}

In this section we give a 
necessary and sufficient condition of finiteness of the set of holes $H$. 
Firstly we will show the necessary and sufficient condition 
% of finiteness of a set of holes $H$ 
in terms of the set of fundamental holes $H_0$.  Then we generalize
the statement, such that it is stated in terms of the minimal Hilbert
basis of $K$.  
Ezra Miller has kindly pointed out to the authors that
many of our results can be proved more succinctly by appropriate
algebraic methods. However for the sake of self-contained presentation
we provide our own proofs and summarize his comments in
Remark \ref{rem:miller1} and Remark \ref{rem:miller2} below.
%, Remark \ref{rem:miller3} ).

First we show that the set of fundamental holes, $H_0$, is finite.
\begin{prop} \label{prop2}
$H_0$ is finite.
\end{prop}
\begin{proof}
Every $\Ba\in \Qsat$ can be written as
\begin{equation}
\label{eq:3}
\Ba = c_1 \Ba_1 + \cdots + c_n \Ba_n,
\end{equation}
where $c_i$'s are nonnegative rational numbers. 
(If $\Ba\in H$, then at least one $c_i$ is not integral.) 
If $c_1 > 1$, then
$\Ba$ can be written as
\[
\Ba = \{(c_1 - \lfloor c_1 \rfloor) \Ba_1 + \cdots + c_n \Ba_n  \} + 
\lfloor c_1 \rfloor \Ba_1 = \tilde \Ba + \lfloor c_1 \rfloor \Ba_1,
\qquad %(\mbox{say})
\]
and $\tilde \Ba = \Ba - \lfloor c_1 \rfloor \Ba_1$.
Therefore
\[
\Qsat \cap (\Ba+(-Q)) \supset \{ \Ba, \tilde \Ba\}
\]
and $\Ba$ is not a fundamental hole.  In this argument we can replace
$c_1$ with any $c_i$, $i\ge 2$.  This shows each fundamental hole has an
expression \eqref{eq:3}, where $0\le c_i \le 1$, $i=1,\ldots,n$.  However %But then
fundamental holes belong to a compact set. Since the lattice points in a
compact set are finite, $H_0$ is finite.
\end{proof}
\begin{rmk}
\label{rem:miller1}
For any field $k$, consider the  semigroup rings $k[Q]$ and $k[\Qsat]$.
Define $M=k[\Qsat]/k[Q]$, which is finitely generated as a
module over $k[Q]$.
$H_0$ is the set of degrees for the minimal generators
of $M$ and therefore $H_0$ is finite.
%  as a module over k[Q].  This follows immediately
%    by definition of H_0.  Since M is a finitely generated
%    k[Q]-module, the set H_0 is always finite.
% $k[Q_{\rm rays}]$, where 
\end{rmk}

Let $H_0 =\{ \By_1, \dots, \By_m\}$. 
Now for each $\By_h \in H_0$ and each $\Ba_i$ define $\bar\lambda_{hi}$ as
follows.  If there exists some $\lambda \in \Z$ such that 
$\By_h + \lambda \Ba_i \in Q$, let
\begin{equation}
\label{eq:lambda-hi}
\bar\lambda_{hi}=
\min \{ \lambda \in \Z  \mid \By_h + \lambda \Ba_i \in Q\}.
\end{equation}
Otherwise define $\bar\lambda_{hi}=\infty$.  Note that
$\bar\lambda_{hi}> 0$  because $\By_h$ is a hole.
Then we have the following result:

\begin{thm}\label{necsuf1}
$H$ is finite if and only if 
$\bar\lambda_{hi} < \infty$ for all $h=1,\dots,m$ and all $i=1,\dots,n$.
\end{thm}

\begin{proof}
For one direction, assume that 
$\bar\lambda_{hi}=\infty$ for some $h$ and $i$.  Then $\By_h + \lambda \Ba_i$,
$\lambda=1,2,\dots$, all belong to $\Qsat$ but do not belong to $Q$.
Therefore they are holes.  Hence $H$ is infinite.

For the other direction, assume that $\bar\lambda_{hi} < \infty$ for all
$h=1,\dots,m$ and all $i=1,\dots,n$.
By \eqref{fundholes}, each hole can be written as
\[
\Bx = \By_h + \sum_{i=1}^n \lambda_i \Ba_i
\]
for some $h$ and $\lambda_i \in \N$, $i=1,\dots,n$.  Now suppose that 
$\lambda_i \ge \bar\lambda_{hi}$ for
some $i$. Then
\[
\By_h + \lambda_i \Ba_i \in Q
\]
and
\[
\Bx = \By_h + \lambda_i \Ba_i + \sum_{j\neq i} \lambda_j \Ba_j \in Q, 
\]
which contradicts that $\Bx$ is a hole.  Therefore if $\Bx$ is a hole,
then $\lambda_i < \bar\lambda_{hi}$ for all $i$. Then
\[
H \subset \{ \By_h + \sum_{i=1}^n \lambda_{hi} \Ba_i \mid h=1,\dots,m, \ 0 \le
\lambda_{hi} < \bar\lambda_{hi}\}.
\]
The right-hand side is finite.
\end{proof}

\begin{rmk}\label{rem:1}
There are several remarks to make.
For each  $1\le i\le n$, let
\[
\tilde Q_{(i)}= \{ \sum_{j\neq i}\lambda_j \Ba_j \mid 
\lambda_j \in \N, \ j\neq i \} 
\]
be the semigroup spanned by $\Ba_j, j\neq i$.  Furthermore write
\[
\bar Q_{(i)} = \Z \Ba_i + \tilde Q_{(i)}.
\]
For each $h$ and $i$,  $\bar\lambda_{hi}$ is finite if and only if 
$
\By_h \in \bar Q_{(i)}.
$
Since $\By_h$ is a hole, actually we only
need to check
\[
\By_h \in (-\N \Ba_i) + \tilde Q_{(i)}.
\]
But $(-\N \Ba_i) + \tilde Q_{(i)}$ is another semigroup, where $\Ba_i$ in $A$ is
replaced by $-\Ba_i$.  Therefore this problem is a standard membership
problem in a semigroup.
\end{rmk}

%Also we only need to check $i$ such that $\Ba_i$ is an extreme ray.
%We have to be somewhat careful in choosing an extreme ray if there are
%multiple vectors in the same direction among the columns of $A$. 
%In this case we take the smallest
%point in $Q$ along the real extreme ray as an extreme ray (although 
%one can take any one of nonzero lattice points in Q along the real ray as 
%an extreme ray).
Also we only need to check $i$ such that $\Ba_i$ is on an extreme ray.  
By a slight abuse of terminology, we simply say that $\Ba_i$ is an extreme ray
if $\Ba_i$ generates an extreme ray of $K$.
If there are multiple columns of $A$ on the same extreme ray, for
definiteness we choose the smallest one, although we can choose any
one of them.
Assume, without loss of generality, that $\{ \Ba_1, \dots,\Ba_k\}$, $k \le
n$, is the set of the extreme rays.  The following corollary says that
we only need to consider $i\le k$.

\begin{cor}\label{necsuf2}
$H$ is finite if and only if 
$\bar\lambda_{hi} < \infty$ for all $h=1,\dots,m$ and all $i=1,\dots,k$.
\end{cor}

\begin{proof}
The first direction is the same as above.  

For the converse direction,
we show that if 
$\bar\lambda_{hi} < \infty$, $1 \le h \le m$, $1 \le i
\le k$, then  $\bar\lambda_{hi} < \infty$, 
$1 \le h \le m$, $k+1 \le i \le n$.   Now any non-extreme ray $\Ba_i$, $i
\ge k+1$, can be written as a nonnegative rational combination of
extreme rays:
\begin{equation}
\label{eq:rays}
\Ba_i = \sum_{j=1}^k  q_{ij} \Ba_j,  \qquad i \ge k+1.
\end{equation}
Let $\bar q_i>0$ denote the l.c.m.\ of the denominators of
$q_{i1},\dots,q_{ik}$.   Then
multiplying both sides by $\bar q_i$, we have
\[
\bar q_i \Ba_i = \sum_{j=1}^k  (\bar q_i q_{ij}) \Ba_j,   \quad \bar q_i q_{ij}
\in \N.
\]
Also note that there is at least one $q_{ij} > 0$, say $q_{ij_0}$.  
Consider $\bar q_i \Ba_i, 2\bar q_i \Ba_i, 3\bar q_i \Ba_i,\dots$. Take 
$\lambda \in \N$ such that
\[
\lambda \bar q_i  q_{ij_0} \ge \bar\lambda_{h j_0}.
\]
Then by \eqref{fundholes}
\[
\By_h + \lambda \bar q_i \Ba_i
= \By_h + \lambda \bar q_i  q_{ij_0} \Ba_{j_0}
+  \sum_{j\neq j_0}^k  \lambda \bar q_i q_{ij} \Ba_j \in Q .
\]
\end{proof}

\begin{rmk}
\label{rem:miller2}
Let $k[Q]$, $k[\Qsat]$ and $M$ be defined as in Remark
\ref{rem:miller1}. 
%$H$ is the set of degrees where $M$ is nonzero.
The number of points in $H$ is the $k$-vector space dimension of $M$.
$H$ is finite if and only if $M$ is Artinian, which proves
Theorem \ref{necsuf1}. Let  
$k[Q_{\rm rays}]$ denote the monoid generated by
the smallest lattice points in $Q$ on the real extreme rays of $Q$.
Then $k[Q]$ is itself finitely generated as a module over the
$k[Q_{\rm rays}]$. This proves Corollary \ref{necsuf2}.
\end{rmk}

Another important point is that we want to state Theorem \ref{necsuf1} in 
terms of Hilbert bases.  Let $B=\{ \Bb_1, \dots, \Bb_L \}$ denote the Hilbert
basis of $K$. %$\Qsat$.  
As above, if $\Bb_l + \lambda \Ba_i \in Q$ for some
$\lambda\in \Z$ let 
\[
\bar \mu_{li} = \min \{ \lambda \in \Z \mid \Bb_l + \lambda \Ba_i
\in Q\}
\]
and $\bar \mu_{li} =\infty$ otherwise. Then we have the following theorem.

\begin{thm}\label{necsuf4}
$H$ is finite if and only if 
$\bar\mu_{li} < \infty$ for all $l=1,\dots,L$ and all $i=1,\dots,n$.
\end{thm}
\begin{proof}
The first direction is the same as the above proofs.

For the converse direction, assume that 
$\bar\mu_{li} < \infty$ for all $l=1,\dots,L$ and all
$i=1,\dots,n$.   Let $\By_h$ be a fundamental hole.
It can be written as a nonnegative integral combination of the
elements of the Hilbert basis
\[
\By_h = \sum_{l=1}^L \alpha_{hl} \Bb_l.
\]
Let
\[
\lambda = \sum_{l=1}^L \alpha_{hl}\bar\mu_{li}.
\]
Then by \eqref{fundholes}
\begin{align*}
\By_h + \lambda \Ba_i &= \sum_{l=1}^L \alpha_{hl} \Bb_l
+ \big(\sum_{l=1}^L \alpha_{hl}\bar\mu_{li}\big) \Ba_i\\
& = \sum_{l=1}^L \alpha_{hl} ( \Bb_l +  \bar\mu_{li} \Ba_i) \in Q.
\end{align*}
This implies $\bar\lambda_{hi} < \infty$ for all $h$ and $i$.
\end{proof}

%Also it is clear that 
As in Corollary \ref{necsuf2}, it is clear that we only need to check extreme rays 
%can be applied to 
%$K$ contained in the Hilbert basis.
among $\Ba_i$'s.

\begin{cor}\label{necsuf3}
$H$ is finite if and only if 
$\bar\mu_{li} < \infty$ for all $l=1,\dots,L$ and all $i=1,\dots,k$.
\end{cor}

\begin{rmk}
In summary,  determining finiteness of $H$ is
straightforward.  We obtain the Hilbert basis $B$ of $\Qsat$. For each
$\Bb\in B\setminus Q$  and for each
extreme $\Ba_i$, we check
\[
\Bb \in  (-\N \Ba_i) + \tilde Q_{(i)}.
\]
\end{rmk}

\begin{ex}
Let $A$ be an integral matrix such that
\[
A = \left(\begin{array}{cccc}
1 & 1 & 1 & 1\\
0 & 1 & 3 & 4\\
\end{array}\right).
\]
%Let $\Ba_i$ be the $i$th column of $A$.
Then $B$ consists of $5$ elements 
\[
B = \{\Bb_1 = (1, 0)^t, \Bb_2 = (1, 1)^t, \Bb_3 = (1, 2)^t, \Bb_4 = (1, 3)^t, 
\Bb_5 = (1, 4)^t\}.
\]
Then we can write $\Bb_3$ as the following:
\begin{eqnarray*}
(1, 2)^t& = & -(1, 0)^t + 2\cdot (1, 1)^t\\
&  = &  (1, 0)^t -(1, 1)^t + (1, 3)^t\\
& = &  (1, 1)^t-(1, 3)^t+(1, 4)^t\\
&  = &  2 \cdot (1, 3)^t - (1, 4)^t.
\end{eqnarray*}
Thus, in this case, we have 
$\bar\mu_{3i} = 1$ for each $i = 1, \, \ldots, 4$
and $\bar\mu_{li} = 0$, where $l \not = 3$ and each 
$i = 1, \, \ldots, 4$.  Thus by Theorem \ref{necsuf4}, the number
of elements in $H$ is finite.
Note that $H$ consists of only one element $\{\Bb_3 = (1, 2)^t\}$.
\end{ex}

\section{Simultaneous finiteness of holes, non-saturation points,  and
minimal saturation points}\label{satpt}

In this section we will show the simultaneous finiteness of holes, 
non-saturation points,  and $S$-minimal saturation points.
As in the previous section let
$\{ \Ba_1, \dots,\Ba_k\}$, $k \le n$, be the set of the extreme rays.
% let $\Ba_1,\dots,\Ba_m$, $m \le n$, be all the extreme rays of $K$.
First, we will show the following lemmas.

\begin{lemma}\label{nonSat}
Suppose that $Q$ is not saturated.
$\Ba\in Q$ is a saturation point if and only if 
$\Ba + \By\in Q$ for all fundamental holes $\By$.
\end{lemma}
\begin{proof}
If $\Ba\in Q $ is a saturation point, then $\Ba + \By \in Q$ for all
$\By\in \Qsat$.  In particular $\Ba + \By\in Q$ for all fundamental
holes $\By$.

Now suppose that $\Ba\in Q$ is not a saturation point.  Then there exists
$\By \in \Qsat$ such that $\Ba + \By$ is a hole.  This $\By$ has to be
a hole, because otherwise $\Ba + \By \in Q$.  $\By$ can be
written as $\By=\By_h + \Bb$ for some fundamental hole $\By_h$ and
$\Bb\in Q$.  Then $\Ba + \By=\Ba + \By_h + \Bb$ and $\Ba+\By_h$
has to be  hole. Therefore we have shown that if $\Ba$ is not a
saturation point, then $\Ba+\By$ is a hole for some fundamental hole $\By$.
\end{proof}

\begin{lemma}\label{equivLemma}
Suppose that $Q$ is  not saturated. Consider any column
$\Ba_i$ of $A$. There 
exists some $n_i \in \N $ such that $n_i \Ba_i \in S$ if and only if
$\bar\lambda_{hi} < \infty$ in (\ref{eq:lambda-hi}) for all
$h=1,\dots,m$.
% , where $\bar\lambda_{hi}$.
\end{lemma}
\begin{proof}
This follows from Lemma \ref{nonSat}.
If $n_i \Ba_i \in S$, $\bar\lambda_{hi} \le
n_i$.  For the other direction take
$n_i = \max_h \bar\lambda_{hi}$. 
\end{proof}

%Now we consider the following three conditions.
Now we consider the following two  conditions.
% To remind a reader that we assume there
% exists $\Bc \in \mathbb{Q}^d$ such that $\Bc \cdot \Ba_i > 0$ for
% $i=1,\ldots,n$, where $\cdot$ is the standard inner product.
\medskip

\noindent
% {\bf Condition 1}\quad   There exists a finite $M>0$ such that every 
% $\Ba \in Q$ with
% $\Bc\cdot \Ba > M$ belongs to $S$.
% %  (i.e.,  $\{ \Ba\in Q \mid \Bc\cdot \Ba > M \}
% % \subset S$).
% \medskip
\noindent
{\bf Condition 1}\quad   For each $\Ba_i$, there exists $n_i>0$ such that $n_i
\Ba_i \in S$.

\medskip
\noindent
{\bf Condition 2}\quad   For each extreme ray $\Ba_i$, $1\le i\le
k$, there exists $n_i>0$ such that $n_i
\Ba_i \in S$.

\medskip

\begin{prop}\label{equiv}
Condition $1$, Condition $2$, and the finiteness of $H$ are 
equivalent.
\end{prop}

\begin{proof}
% Suppose that Condition $1$ holds. For each $\Ba_i$, let  $n_i > M/(\Bc\cdot
% \Ba_i)$, then  $ \Bc \cdot (n_i \Ba_i) > M$ and  $n_i \Ba_i \in S$.
% Therefore Condition $2$ holds.  
% Condition $3$ trivially holds if Condition
% $2$ holds.
Condition $1$ trivially implies Condition $2$.
On the other hand suppose that Condition $2$ holds.  Then each
non-extreme $\Ba_i$, $k < i \le n$, can be written as 
(\ref{eq:rays}).
% a non-negative
% rational combination of extreme rays:
% \begin{equation}
% \label{eq:rational}
% \Ba_j = \sum_{i=1}^m \lambda_i \Ba_i.
% \end{equation}
As above let $\bar q_i>0$ denote the l.c.m.\ of the denominators of
$q_{i1},\dots,q_{ik}$  and let
$n_i = \bar q_i \times n_1 \times  \dots \times n_k$, then
$n_i \Ba_i \in S$ and Condition 1 holds.
\comment{
Next suppose that Condition $2$ holds.  Then every
non-saturation point $\Ba\in \bar S$ has an expression
\[
\Ba = \lambda_1 \Ba_1 + \dots + \lambda_n \Ba_n, \qquad
  0\le \lambda_i <  n_i, \forall i .
\]
Therefore $\bar S$ is a subset of a compact set and hence finite.  Choose
$M$ such that
\[
M > \max_{\Ba \in \bar S} \Bc \cdot \Ba.
\]
Then Condition $1$ holds.}

%Finally 
Now we show the equivalence between the finiteness of $H$ and
the other two conditions.  Using Lemma \ref{equivLemma}, Condition 1 is 
equivalent to the condition in Theorem \ref{necsuf1}.  Also Condition 2 is 
equivalent to the condition in Corollary \ref{necsuf2}.
\end{proof}

Now we prove Theorem \ref{condition}. 
%In the theorem $\cone(S)$ 
%denotes the set of finite nonnegative real combinations of
%elements of $S$ and ``rational polyhedral cone'' is a closed cone defined
%by rational linear weak inequalities.

\comment{
\begin{thm}\label{condition}
%If Condition $1$ holds, then $H, \bar S$ and
%$\min(S;S)$ are all finite.  If Condition $1$ does not hold, then $H, \bar
%S$ and $\min(S;S)$ are all infinite.
%$H, \bar S$, and $\min(S;S)$ are simultaneously finite or simultaneously
%infinite.
The following statements are equivalent.  \begin{enumerate}
\item\label{equiv1} $\min(S;S)$ is finite. 
\item\label{equiv2} $\cone(S)$ is a rational polyhedral cone.
\item\label{equiv3} There is some $s\in S$ on every extreme ray of $K$.
\item\label{equiv4} $H$ is finite.
\item\label{equiv5} $\bar S$ is finite.
\end{enumerate}
\end{thm}
}

\begin{proof}[Proof for Theorem \ref{condition}]
\noindent
\ref{equiv1}. $\Longleftrightarrow$ \ref{equiv2}. :  
% $\min(S;S)$ is finite or equivalently, the monoid $S$, spanned by $\min(S;S)$,
% has a finite integral generating set.
% We then apply Theorem 1.1 (b)
% in [Hemmecke and Weismantel, 2006] (also we can apply Theorem 4 in [Jeroslow, 1978], which 
% is a corollary of Theorem 1.1 (b) in [Hemmecke and Weismantel, 2006]).
$\min(S;S)$ is an integral generating set of the monoid $S \cup \{0\}$. 
We then apply Theorem 1.1 (b) of 
[\cite{hemmeke-weismantel}] %[Hemmecke and Weismantel, 2006] 
or Theorem 4 in [\cite{jeroslow1978}]. % [Jeroslow, 1978].

\smallskip
\noindent
\ref{equiv2}. $\Longleftrightarrow$ \ref{equiv3}. : 
If $\cone(S)$ is not polyhedral, there must be an extreme ray $e$ of $K$ not in
$\cone(S)$, since $K$ is polyhedral. Thus, $e\cap S=\emptyset$.

If $\cone(S)$ is polyhedral, then it is a 
rational polyhedron and has a finite integral generating set. Thus, 
by Theorem 8.8 in [\cite{bertsimas-weismantel}] %[Bertsimas and Weismantel, 2005], 
the polyhedron $\cone(S)$ contains all lattice points from its recession cone 
$K$
%but finitely many lattice
%points from its recession cone $K$,
and $(K \setminus \cone(S))\cap \Z^d$ is finite, which in this case can only happen if
$\cone(S)=K$. Thus, there is a point from $S$ on each extreme ray of $K$.
\smallskip

\noindent
\ref{equiv3}. $\Longleftrightarrow$ \ref{equiv4}. : The statement \ref{equiv3}.
is equivalent to Condition $2$.  Thus, the proof follows
directly by Proposition \ref{equiv}.
\smallskip

\noindent
\ref{equiv4}. $\Longleftrightarrow$ \ref{equiv5}. :
Suppose that $H$ is finite.  Then
% Condition $1$ holds and as in the proof of Proposition
% \ref{equiv}
% $\bar S$ is finite.
by Condition 1, it is easy to see that $\bar S$ is contained in a
compact set and hence $\bar S$ is finite.
For the opposite implication, suppose that $H$ is infinite.
Since Condition 1 does not hold, 
there exists some $i$ such that $n \Ba_i \not\in S$ for all $n\in \N$.
Then $\{ \Ba_i, 2\Ba_i,3\Ba_i,\dots \} \subset \bar S$ and $\bar S$ is
infinite.
\end{proof}

%Now we will prove Theorem \ref{mainthm}.

%\begin{proof}
%By Proposition \ref{prop0}, we know that if $H$ is finite then $\min(S; S)$ is
%finite.  Now we will prove the other direction.
%We argue by contradiction.  Suppose $\min(S; S)$ is finite but
%$H$ is infinite.
%\end{proof}

%%%%%%%%%%%%%%%%%

Now we consider the generators 
%of $\min(S;Q_{\rm sat})$.  and  
$\min(S;Q)$ and
we prove that $\min(S;Q)$ is always finite.
% $\min(S;Q_{\rm sat})$ is finite.
Then by (\ref{eq:2})  $\min(S;\Qsat)$ is always finite as well.
Note that 
the multi-dimensional Frobenius problem can be stated as computing
the sets   $\min(S;Q)$ and $\min(S;Q_{\rm sat})$.

% \begin{prop} \label{prop1}
% $\min(S;Q_{\rm sat})$ is finite.
% \end{prop}
% \begin{proof}
% We can assume without loss of generality that $\Z^d = L$.
% Thus $Q_{\rm sat} = K \cap \Z^d$.  Theorem 16.4 in \cite{Schrijver1986} 
% states that $K$ has a unique minimal
% Hilbert basis (under our assumption of pointedness of $Q$).
% Thus, $Q_{\rm sat}$ is a finitely generated monoid.
% Let $\{{\rm h}_1, \cdots, {\rm h}_q\}$, where $q \leq d$ is a positive 
% integer, be the Hilbert basis. 
% Consider the algebra, $k[Q_{\rm sat}] := k[t^{{\rm h}_1}, \cdots,
% t^{{\rm h}_q}]$, where $k$ is any algebraic field.  
% Then $k[Q_{\rm sat}]$ is a finitely generated $k$-algebra and therefore
% a Noetherian ring by a corollary of Hilbert's basis theorem (Corollary 1.3) in 
% \cite{Eisenbud1995}.
% By the definition of $S$, one notice that $a + b \in S, \ \forall a \in S, \ 
% \forall b \in Q_{\rm sat}$.  So $I_S := <t^\beta: \beta \in S>$ is an ideal
% in $k[Q_{\rm sat}]$.  By the definition of a Noetherian ring,
% every ideal is generated finitely.  Thus $I_S$ is generated finitely.
% \end{proof}

\begin{prop} 
$\min(S;Q)$ is finite.
\end{prop}
\begin{proof}
Note that $Q$ is a finitely generated monoid. 
Consider the algebra, $k[Q] := k[t^{\Ba_1}, \cdots,
t^{\Ba_n}]$, where $k$ is any algebraic field.
Then $k[Q]$ is a finitely generated $k$-algebra by Proposition 2.5 in 
[\cite{bruns-gubeladze}] %[Bruns and Gubeladze, 2006] 
and therefore
a Noetherian ring by a corollary of Hilbert's basis theorem 
(Corollary 1.3 in %[Eisenbud, 1995]
[\cite{Eisenbud1995}]).  Since $I_S:= <t^\beta: \beta \in S>$ is an ideal in 
$k[Q]$, we are done.
\end{proof}
A combinatorial proof of this proposition is given in [\cite{Hemmecke}].
%Finally we point out that the difference  of $\min (S;Q)$ and 
%$\min(S;\Qsat)$ is contained in the set of fundamental holes in the
%following sense.

\begin{prop}
\begin{equation}
\min (S;Q)\subset \min(S;\Qsat) + (H_0\cup \{0\}).
\end{equation}
\end{prop}

\begin{proof}
%Recall that $H_0$ contains the origin $0$.  <- change of definition
Let $\Ba\in \min(S;Q)$.  We want to show that
$\Ba$ can be written as $\Ba = \tilde\Ba + \Bb$, where $\tilde \Ba \in
\min(S;\Qsat)$ and $\Bb\in H_0\cup \{0\}$. If $\Ba$ itself belongs to
$\min(S;\Qsat)$, then take $\Ba=\tilde \Ba$ and $\Bb=0$.
Otherwise, if $\Ba \not\in \min(S;\Qsat)$, then
by definition of $\Qsat$-minimality there  exists $\Ba' \in S$ such that
$0\neq \Ba - \Ba'\in \Qsat$.  If $\Ba' \not\in \min(S;\Qsat)$, then we can do
the same operation to $\Ba'$.  This operation has to stop in finite steps
and we arrive at $\tilde \Ba\in \min(S;\Qsat)$ such that
$\Bb=\Ba - \tilde \Ba \in \Qsat$.  If this $\Bb\not\in H_0$, then
there exists $\Bc \in Q$, $\Bc\neq 0$, such that
$\Bb - \Bc\in \Qsat$.  Then
\[
\Ba = \tilde \Ba + \Bb = \tilde \Ba + (\Bb - \Bc) + \Bc,
\]
where $\tilde \Ba \in S$, $\Bb- \Bc \in \Qsat$.  
%By monotonicity of $S$,
Since $S + \Qsat \subset S$,
$\tilde \Ba + (\Bb - \Bc)\in S$.  But this contradicts $\Ba\in \min(S;Q)$.
\end{proof}

\section{Applications to contingency tables}\label{sec:contingency-table}
An {\bf $s$-way contingency table of size $n_1\times \cdots \times
  n_s$}  is an array of nonnegative
integers $v=(v_{i_1,\dots,i_s})$, $1\leq i_j\leq n_j$. For $0\leq r
<s$, an {\bf $r$-marginal} of $v$ is any of the $s\choose r$ possible
$r$-way tables obtained by summing the entries over all but $r$
indices. 
In this section we apply our theorem to some examples including 
$2\times 2\times 2\times 2$ tables with $2$-marginals 
% (the complete  graph with $4$ nodes and  with levels of $2$ on each node, $K4$), 
and
 $2\times 2\times 2\times 2$ tables with three $2$-marginals and a $3$-marginal
($[12][13][14][234]$).  Also we apply our theorem to
three-way contingency tables from %[Vlach, 1986].
%\cite{Irving1994}.
[\cite{Vlach1986}].
To compute minimal Hilbert bases of cones, we used 
{\tt normaliz} [\cite{brunskoch}]  %[Bruns and Koch, 2001] 
and to compute each hyperplane representation and 
vertex representation we used 
{\tt CDD} [\cite{fukuda}] % [Fukuda, 2005] 
and 
{\tt lrs} [\cite{avis}]. %[Avis, 2005].
Also we used {\tt 4ti2} %\cite{hemmecke2005} [Hemmecke et al, 2005] 
[\cite{Hemmecke+Hemmecke+Malkin:2005}]
to compute matrix $A$ for the system.

\subsection*{$2\times 2\times 2\times 2$ tables}
\subsubsection*{$2\times 2\times 2\times 2$ tables with $2$-marginals}
First, we would like to show some simulation results with 
$2\times 2\times 2\times 2$ tables with $2$-marginals, which can be seen as 
the complete graph with $4$ nodes $K4$ and  with $2$ states on each node.  
The semigroup of $K4$ has $16$ generators
$\Ba_1,\dots,\Ba_{16}$ in $\Z^{24}$ (without removing redundant rows) such that
% \begin{verbatim}
% 1 0 0 0 1 0 0 0 1 0 0 0 1 0 0 0 1 0 0 0 1 0 0 0
% 0 1 0 0 0 1 0 0 0 1 0 0 1 0 0 0 1 0 0 0 1 0 0 0
% 0 0 1 0 1 0 0 0 1 0 0 0 0 1 0 0 0 1 0 0 1 0 0 0
% 0 0 0 1 0 1 0 0 0 1 0 0 0 1 0 0 0 1 0 0 1 0 0 0
% 1 0 0 0 0 0 1 0 1 0 0 0 0 0 1 0 1 0 0 0 0 1 0 0
% 0 1 0 0 0 0 0 1 0 1 0 0 0 0 1 0 1 0 0 0 0 1 0 0
% 0 0 1 0 0 0 1 0 1 0 0 0 0 0 0 1 0 1 0 0 0 1 0 0
% 0 0 0 1 0 0 0 1 0 1 0 0 0 0 0 1 0 1 0 0 0 1 0 0
% 1 0 0 0 1 0 0 0 0 0 1 0 1 0 0 0 0 0 1 0 0 0 1 0
% 0 1 0 0 0 1 0 0 0 0 0 1 1 0 0 0 0 0 1 0 0 0 1 0
% 0 0 1 0 1 0 0 0 0 0 1 0 0 1 0 0 0 0 0 1 0 0 1 0
% 0 0 0 1 0 1 0 0 0 0 0 1 0 1 0 0 0 0 0 1 0 0 1 0
% 1 0 0 0 0 0 1 0 0 0 1 0 0 0 1 0 0 0 1 0 0 0 0 1
% 0 1 0 0 0 0 0 1 0 0 0 1 0 0 1 0 0 0 1 0 0 0 0 1
% 0 0 1 0 0 0 1 0 0 0 1 0 0 0 0 1 0 0 0 1 0 0 0 1
% 0 0 0 1 0 0 0 1 0 0 0 1 0 0 0 1 0 0 0 1 0 0 0 1
% \end{verbatim}
\begin{verbatim}
 1 0 0 0 1 0 0 0 1 0 0 0 1 0 0 0
 0 1 0 0 0 1 0 0 0 1 0 0 0 1 0 0
 0 0 1 0 0 0 1 0 0 0 1 0 0 0 1 0
 0 0 0 1 0 0 0 1 0 0 0 1 0 0 0 1
 1 0 1 0 0 0 0 0 1 0 1 0 0 0 0 0
 0 1 0 1 0 0 0 0 0 1 0 1 0 0 0 0
 0 0 0 0 1 0 1 0 0 0 0 0 1 0 1 0
 0 0 0 0 0 1 0 1 0 0 0 0 0 1 0 1
 1 0 1 0 1 0 1 0 0 0 0 0 0 0 0 0
 0 1 0 1 0 1 0 1 0 0 0 0 0 0 0 0
 0 0 0 0 0 0 0 0 1 0 1 0 1 0 1 0
 0 0 0 0 0 0 0 0 0 1 0 1 0 1 0 1
 1 1 0 0 0 0 0 0 1 1 0 0 0 0 0 0
 0 0 1 1 0 0 0 0 0 0 1 1 0 0 0 0
 0 0 0 0 1 1 0 0 0 0 0 0 1 1 0 0
 0 0 0 0 0 0 1 1 0 0 0 0 0 0 1 1
 1 1 0 0 1 1 0 0 0 0 0 0 0 0 0 0
 0 0 1 1 0 0 1 1 0 0 0 0 0 0 0 0
 0 0 0 0 0 0 0 0 1 1 0 0 1 1 0 0
 0 0 0 0 0 0 0 0 0 0 1 1 0 0 1 1
 1 1 1 1 0 0 0 0 0 0 0 0 0 0 0 0
 0 0 0 0 1 1 1 1 0 0 0 0 0 0 0 0
 0 0 0 0 0 0 0 0 1 1 1 1 0 0 0 0
 0 0 0 0 0 0 0 0 0 0 0 0 1 1 1 1
\end{verbatim}

Remember that the columns of the given array are the generators of the
semigroup.
All of these vectors are extreme rays of the cone, which we verified via 
{\tt cddlib} [\cite{fukuda}]. %[Fukuda, 2005].
The Hilbert basis of the cone generated by these $16$ vectors contains $17$
vectors $\Bb_1,\dots,\Bb_{17}$. The first 16 vectors are the same as
$\Ba_i$, i.e.\ $\Bb_i=\Ba_i$, $i=1,\dots,16$. The 17-th vector $\Bb_{17}$
is 
\[
\Bb_{17}=(1\ 1\ \dots \ 1)^t\]
consisting of all 1's.
Thus, $\Bb_{17} \not \in Q$.  
Then we set the $16$ systems of linear equations such that:
\begin{eqnarray*}
P_j: & \Bb_1 x_1 + \Bb_2 x_2 + \cdots + \Bb_{16} x_{16} = \Bb_{17}\\
 &x_j \in \Z_-, \, \, x_i \in  \Z_+, \, \mbox{ for } i \not = j,
\end{eqnarray*}
for $j = 1, 2, \cdots , 16$.
We solved these systems via {\tt lrs} and {\tt LattE}
[\cite{lattemanual}].
% [De Loera et al., 2003].  
Then we have:
{\allowdisplaybreaks
\begin{eqnarray*}
\Bb_{17} & = & -\Bb_1 + \Bb_2 + \Bb_3 + \Bb_5 + \Bb_9 + \Bb_{16},\\
\Bb_{17} & = & \Bb_1 - \Bb_2 + \Bb_4 + \Bb_6 + \Bb_{10} + \Bb_{15},\\
\Bb_{17} & = & \Bb_1 - \Bb_3 + \Bb_4 + \Bb_7 + \Bb_{11} + \Bb_{14},\\
\Bb_{17} & = & \Bb_2 + \Bb_3 - \Bb_4 + \Bb_8 + \Bb_{12} + \Bb_{13},\\
\Bb_{17} & = & \Bb_1 - \Bb_5 + \Bb_6 + \Bb_7 + \Bb_{12} + \Bb_{13},\\
\Bb_{17} & = & \Bb_2 + \Bb_5 - \Bb_6 + \Bb_8 + \Bb_{11} + \Bb_{14},\\
\Bb_{17} & = & \Bb_3 + \Bb_5 - \Bb_7 + \Bb_8 + \Bb_{10} + \Bb_{15},\\
\Bb_{17} & = & \Bb_4 + \Bb_6 + \Bb_7 - \Bb_8 + \Bb_{9} + \Bb_{16},\\
\Bb_{17} & = & \Bb_1 + \Bb_8 - \Bb_9 + \Bb_{10} + \Bb_{11} + \Bb_{13},\\
\Bb_{17} & = & \Bb_2 + \Bb_7 + \Bb_9 - \Bb_{10} + \Bb_{12} + \Bb_{14},\\
\Bb_{17} & = & \Bb_3 + \Bb_6 + \Bb_9 - \Bb_{11} + \Bb_{12} + \Bb_{15},\\
\Bb_{17} & = & \Bb_4 + \Bb_5 + \Bb_{10} + \Bb_{11} - \Bb_{12} + \Bb_{16},\\
\Bb_{17} & = & \Bb_4 + \Bb_5 + \Bb_{9} - \Bb_{13} + \Bb_{14} + \Bb_{15},\\
\Bb_{17} & = & \Bb_3 + \Bb_6 + \Bb_{10} + \Bb_{13} - \Bb_{14} + \Bb_{16},\\
\Bb_{17} & = & \Bb_2 + \Bb_7 + \Bb_{11} + \Bb_{13} - \Bb_{15} + \Bb_{16},\\
\Bb_{17} & = & \Bb_1 + \Bb_8 + \Bb_{12} + \Bb_{14} + \Bb_{15} - \Bb_{16}.
\end{eqnarray*}
}
Thus by Theorem \ref{necsuf4}, the number of elements in $H$ is finite.

\subsubsection*{$2\times 2\times 2\times 2$ tables with $2$-marginals and a $3$-marginal}

Now we consider $2 \times 2 \times 2 \times 2$ tables with three 
$2$-marginals and a $3$-marginal as the simplicial complex on $4$ nodes 
$[12][13][14][234]$ and with $2$ states on each node.

%\begin{figure}[!htb]
%  \centering
%  \setlength{\unitlength}{0.08mm}
%  \begin{picture}(800,800)
%    \thicklines
%    % x-axis
%    \put(20,20){\line(1, 0){760}} \put(0,0){\makebox(0,0){1}}
%    \put(20,20){\line(0, 1){750}} \put(800,0){\makebox(0,0){2}}
%    \put(780,780){\line(-1, -1){750}} \put(20,780){\line(1, -1){750}} 
%    \put(0,800){\makebox(0,0){4}} \put(780,780){\line(-1, 0){765}}
%    \put(780,780){\line(0, -1){765}} \put(800,800){\makebox(0,0){3}}
%  \end{picture}
%\caption{The simplicial complex on 4 nodes}\label{fig1}
%\end{figure}

After removing redundant rows (using {\tt cddlib}),
$2 \times 2 \times 2 \times 2$ tables with $2$-marginals and a
$3$-marginal has the $12 \times 16$  matrix $A$.  Thus the semigroup is 
generated by $16$ vectors in $\Z^{12}$ such that:
% \begin{verbatim}
% 1 0 0 0 1 0 1 0 1 0 0 0 
% 0 1 0 0 0 1 0 1 1 0 0 0 
% 0 0 1 0 1 0 1 0 0 1 0 0 
% 0 0 0 1 0 1 0 1 0 1 0 0 
% 1 0 0 0 0 0 1 0 0 0 1 0 
% 0 1 0 0 0 0 0 1 0 0 1 0 
% 0 0 1 0 0 0 1 0 0 0 0 0 
% 0 0 0 1 0 0 0 1 0 0 0 0 
% 1 0 0 0 1 0 0 0 0 0 0 1 
% 0 1 0 0 0 1 0 0 0 0 0 1 
% 0 0 1 0 1 0 0 0 0 0 0 0 
% 0 0 0 1 0 1 0 0 0 0 0 0 
% 1 0 0 0 0 0 0 0 0 0 0 0 
% 0 1 0 0 0 0 0 0 0 0 0 0 
% 0 0 1 0 0 0 0 0 0 0 0 0 
% 0 0 0 1 0 0 0 0 0 0 0 0 
% \end{verbatim}
\begin{verbatim}
 1 0 0 0 1 0 0 0 1 0 0 0 1 0 0 0
 0 1 0 0 0 1 0 0 0 1 0 0 0 1 0 0
 0 0 1 0 0 0 1 0 0 0 1 0 0 0 1 0
 0 0 0 1 0 0 0 1 0 0 0 1 0 0 0 1
 1 0 1 0 0 0 0 0 1 0 1 0 0 0 0 0
 0 1 0 1 0 0 0 0 0 1 0 1 0 0 0 0
 1 0 1 0 1 0 1 0 0 0 0 0 0 0 0 0
 0 1 0 1 0 1 0 1 0 0 0 0 0 0 0 0
 1 1 0 0 0 0 0 0 0 0 0 0 0 0 0 0
 0 0 1 1 0 0 0 0 0 0 0 0 0 0 0 0
 0 0 0 0 1 1 0 0 0 0 0 0 0 0 0 0
 0 0 0 0 0 0 0 0 1 1 0 0 0 0 0 0
\end{verbatim}

% Remember that the columns of the given array are the generators of the
% semigroup.
%Note that 
All of these vectors are extreme rays of the cone (verified via 
{\tt cddlib}).
The Hilbert basis of the cone generated by these $16$ vectors
consists of these 16 vectors and two additional vectors
% vectors as the following:
% \begin{verbatim}
% 1 0 0 0 1 0 1 0 1 0 0 0
% 0 1 0 0 0 1 0 1 1 0 0 0
% 0 0 1 0 1 0 1 0 0 1 0 0
% 0 0 0 1 0 1 0 1 0 1 0 0
% 1 0 0 0 0 0 1 0 0 0 1 0
% 0 1 0 0 0 0 0 1 0 0 1 0
% 0 0 1 0 0 0 1 0 0 0 0 0
% 0 0 0 1 0 0 0 1 0 0 0 0
% 1 0 0 0 1 0 0 0 0 0 0 1
% 0 1 0 0 0 1 0 0 0 0 0 1
% 0 0 1 0 1 0 0 0 0 0 0 0
% 0 0 0 1 0 1 0 0 0 0 0 0
% 1 0 0 0 0 0 0 0 0 0 0 0
% 0 1 0 0 0 0 0 0 0 0 0 0
% 0 0 1 0 0 0 0 0 0 0 0 0
% 0 0 0 1 0 0 0 0 0 0 0 0
% 1 1 1 1 1 1 1 1 1 0 0 0
% 1 1 1 1 1 1 1 1 0 1 1 1
% \end{verbatim}
\[
\Bb_{17}=(1\ 1\ 1\ 1\ 1\  1\  1\  1\  1\  0\  0\  0)^t, \qquad
\Bb_{18}=(1\  1\  1\  1\  1\  1\  1\  1\  0\  1\  1\  1)^t.
\]
% Let $\Bb_i$ be the $i$th elements in the Hilbert basis above (i.e. the $i$th
% row of the list above).  
Thus, $\Bb_{17}, \, \Bb_{18} \not \in Q$.  
Then we set the system of linear equations such that:
\begin{eqnarray*}
\Bb_1 x_1 + \Bb_2 x_2 + \cdots + \Bb_{16} x_{16} = \Bb_{17}\\
x_1 \in \Z_-, \, \, x_i \in  \Z_+, \, \mbox{ for } i = 2, \cdots , 16.
\end{eqnarray*}
We solved the system via {\tt lrs} and {\tt CDD}. 
We noticed that this system has no real solution 
(infeasible).  This means that 
\[
\Bb_{17} \not \in  (-\N \Ba_1) + \tilde Q_{(1)}.
\]
Thus by Theorem \ref{necsuf4}, the number of elements in $H$ is infinite.

%\subsection*{Simulations with three dimensional tables}
\subsection*{Results on  three-way tables}

% In this section we study the semigroups of defining matrices for some 
% three dimensional tables.  For $2 \times J \times K$ tables with
% $2$-marginals for $J, \, K \in \N$, 
% the semigroup of each defining matrix is saturated
% since $2 \times J \times K$ is the Lawrence lifting of $J \times K$, 
% which is also unimodular (and Lawrence lifting of unimodular is unimodular
% \cite{Sturmfels1996}).
% %, \cite{Sturmfels1998}% [Sturmfels, 1996, Sturmfels, 1998]).  
% Thus, an example in
% \cite{Irving1994} %[Irving, 1994]
% does not have a table with integral entries as well as a table with real 
% entries because this is unimodular.

% For $3 \times 3 \times J$ tables with
% $2$-marginals for $J \in \N$, 
% we do not know that the semigroup of 
% each defining matrix is saturated for all $J \in \N$.
Results on the saturation of 3-DIPTP are summarized in Theorem 6.4 of 
[\cite{ohsugi-hibi-contingency-tables-2006}]. 
They show that a normality (i.e., $Q$ is saturated) or non-normality
(i.e., $Q$ is not saturated) 
of $Q$ is not known only for the following three cases:
\[
5\times 5\times 3, \quad 5\times 4 \times 3, \quad 4\times 4\times 3.
\]
All $2 \times J \times K$ tables are unimodular and hence saturated.
This means that there is no hole in $Q$, and
thus a $2\times 2\times 2$ example in [\cite{Irving1994}] %[Irving, 1994]
is not a hole.
% does not have a table with integral entries as well as a table with real 
% entries because this is unimodular.
All $3 \times 3 \times J$ tables are saturated by the result of 
\cite{sullivant2004}.

\begin{figure}[h]
\begin{center}
     \includegraphics[scale=.4]{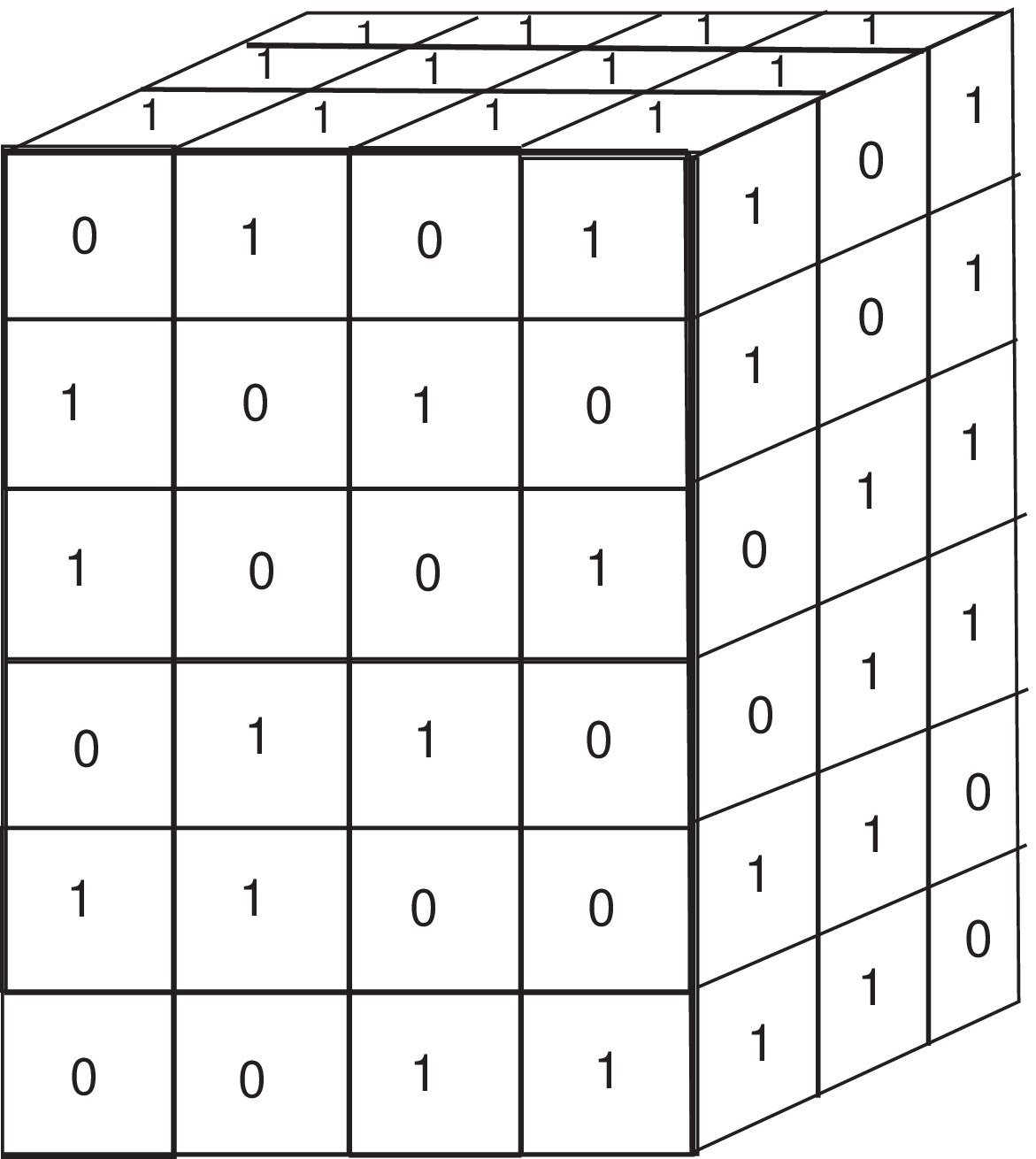}
 \end{center}
\caption{An example of $3 \times 4 \times 6$ table such that the given marginal
condition is a hole of the semigroup.}\label{noninteger}
\end{figure}

For $3 \times 4 \times 6$ tables with $2$-marginals, \cite{Vlach1986} showed an
example which has a table with nonnegative real entries, but does not
have a table with nonnegative integer entries. 
%[\cite{Vlach1986}]. %[Vlach, 1986].
This example can be found in Figure \ref{noninteger}.  Actually it is
a particular example of Lemma 6.1 of
[\cite{ohsugi-hibi-contingency-tables-2006}].
\cite{aoki-takemura-metr03-38} presents a non-squarefree indispensable
move $\Bz=\Bz^+ - \Bz^-$ of size $3 \times 4 \times 6$, where $2$
appears both in the positive part $\Bz^+$ and the negative part $\Bz^-$.
For this $\Bz$ there exist two standard coordinate vectors $\Be_1,
\Be_2$ such that 
\[
\Bu = \Bz^+ - 2\Be_1 \ge 0, \quad 
\Bv = \Bz^- - 2\Be_2 \ge 0.
\]
In this case Lemma 6.1 of [\cite{ohsugi-hibi-contingency-tables-2006}]
proves that $\Bb = A (\Bu + \Bv)/2 \in \Qsat$ is a hole and this
corresponds to Vlach's example.

Using Vlach's example, one can also show that
$3 \times 4 \times 7$ tables and bigger tables have infinitely many 
holes.  We take the example in Figure \ref{noninteger}.  Then we
embed the table in a $3 \times 4 \times 7$ table.  
Then we put a single arbitrary positive integer $c$ at just one place of the 
seventh $3\times 4$ slice.  This positive integer is uniquely determined by
2-marginals of the seventh slice alone % (Table \ref{noninteger2}).  
(Table 1).
Thus for each choice of $c$
the beginning $3 \times 4 \times 6$ part remains  to
be a hole.  Since $c$ is arbitrary, 
$3 \times 4 \times 7$ table has infinite number of holes.
% take any integer in the cell of the last corner of 
% the  $4 \times 5 \times 7$ table otherwise all $0$.  
% this family of
% tables are in the set of holes. 

\begin{table}[h]
\begin{center}
%\caption{Averages of computational experiences} 
\begin{tabular}{|c|c|c|c|c||c|} \hline
& & & & & sum \\\hline
& $c$ & 0 & 0 & 0 & $c$\\\hline
& 0  & 0 & 0 & 0 & 0\\\hline
& 0  & 0 & 0 & 0 & 0\\\hline \hline
sum & $c$ & 0 & 0 & 0 & $c$\\\hline
\end{tabular}
\caption{the $7$-th $3 \times 4$ slice is uniquely determined by its row and 
its column sums.  $c$ is an arbitrary positive integer.}
\end{center}\label{noninteger2}
\end{table}
We can generalize this idea as follows.
Let $A_1$ denote the integer matrix corresponding to problem of a
smaller size.  Suppose that $A$ for a larger problem can be written as
a partitioned matrix 
\[
\begin{pmatrix}
 A_1 & 0 \\
 0   & A_2 \\
 A_3 & A_4 \\
\end{pmatrix},
\]
where $A_3$ and $A_4$ are arbitrary.  We consider the case that for
$A_1$ there exists a hole.  Now consider the semigroup associated with
$A_2$.  We assume that there exists infinite number of one-element
fibers for the semigroup associated with $A_2$.   This is usually the
case, because the fibers on the extreme ray for $A_2$ is all
one-element fibers, under the condition that $A_2$ does not contain
%multiple 
more than one extreme rays in the same direction.

Under these assumptions consider the equation
\[
\begin{pmatrix} t_1 \\ t_2 \\ t_3 \end{pmatrix}
= 
\begin{pmatrix}
 A_1 & 0 \\
 0   & A_2 \\
 A_3 & A_4 \\
\end{pmatrix} 
\begin{pmatrix} x_1 \\ x_2
\end{pmatrix},
\]
where $t_1$ is a hole for $A_1$, $t_2$ is any of the one-element
fibers for $A_2$ and $t_3$ is chosen to satisfy the equation. Then
$(t_1, t_2, t_3)^t$ is a hole for each $t_2$.  Therefore there
exist infinite number of holes for the larger problem.
\begin{ex}
Let $A_1$ be an integral matrix such that
\[
A_1 = \left(\begin{array}{cccc}
1 & 1 & 1 & 1\\
0 & 1 & 3 & 4\\
\end{array}\right).
\]
and let $A_2 = (1)$. From Example \ref{example1},
$H$ consists of only one element $\{t_1 = (1, 2)^t\}$ and
with $A_2$ we can find a family of infinite number of one-element fibers, 
namely $F_{c} := \{c\}$, where $c$ is an arbitrary positive integer.  
Let $t_2 = c$. Then we have a matrix $A$ such that:
\[
A = \left(\begin{array}{cc} A_1 & 0 \\
0 & A_2\\
\end{array}\right)
= \left(\begin{array}{ccccc}
1 & 1 & 1 & 1 & 0\\
0 & 1 & 3 & 4 & 0\\
0 & 0 & 0 & 0 & 1\\
\end{array}\right).
\]
Note that $(t_1, t_2)^t = (1, 2, c)^t$ is a hole for each $t_2 = c$.
Thus, since $c$ is an arbitrary positive integer,
there exist infinitely many holes for the semigroup 
generated by the columns of the matrix $A$.
\end{ex}

\comment{
\subsubsection{$3 \times 3 \times J$ tables with $2$-marginals}
In this section we investigate the normality of the semigroup of 
$3 \times 3 \times J$ tables for $J \leq 5$.  
For example, $3 \times 3 \times 3$ with $2$-marginals, 
after removing redundant rows (we removed redundant rows using {\tt cddlib}),
the constraint matrix becomes a $19 \times 27$ matrix.  Thus the semigroup
is generated by $27$ vectors in $\Z^{19}$ such that:
% 1 0 0 0 0 0 0 0 0 1 0 0 0 0 0 1 0 0 0 
% 0 1 0 0 0 0 0 0 0 0 1 0 0 0 0 1 0 0 0 
% 0 0 1 0 0 0 0 0 0 0 0 1 0 0 0 1 0 0 0 
% 0 0 0 1 0 0 0 0 0 1 0 0 0 0 0 0 1 0 0 
% 0 0 0 0 1 0 0 0 0 0 1 0 0 0 0 0 1 0 0 
% 0 0 0 0 0 1 0 0 0 0 0 1 0 0 0 0 1 0 0 
% 0 0 0 0 0 0 1 0 0 1 0 0 0 0 0 0 0 0 0 
% 0 0 0 0 0 0 0 1 0 0 1 0 0 0 0 0 0 0 0 
% 0 0 0 0 0 0 0 0 1 0 0 1 0 0 0 0 0 0 0 
% 1 0 0 0 0 0 0 0 0 0 0 0 1 0 0 0 0 1 0 
% 0 1 0 0 0 0 0 0 0 0 0 0 0 1 0 0 0 1 0 
% 0 0 1 0 0 0 0 0 0 0 0 0 0 0 1 0 0 1 0 
% 0 0 0 1 0 0 0 0 0 0 0 0 1 0 0 0 0 0 1 
% 0 0 0 0 1 0 0 0 0 0 0 0 0 1 0 0 0 0 1 
% 0 0 0 0 0 1 0 0 0 0 0 0 0 0 1 0 0 0 1 
% 0 0 0 0 0 0 1 0 0 0 0 0 1 0 0 0 0 0 0 
% 0 0 0 0 0 0 0 1 0 0 0 0 0 1 0 0 0 0 0 
% 0 0 0 0 0 0 0 0 1 0 0 0 0 0 1 0 0 0 0 
% 1 0 0 0 0 0 0 0 0 0 0 0 0 0 0 0 0 0 0 
% 0 1 0 0 0 0 0 0 0 0 0 0 0 0 0 0 0 0 0 
% 0 0 1 0 0 0 0 0 0 0 0 0 0 0 0 0 0 0 0 
% 0 0 0 1 0 0 0 0 0 0 0 0 0 0 0 0 0 0 0 
% 0 0 0 0 1 0 0 0 0 0 0 0 0 0 0 0 0 0 0 
% 0 0 0 0 0 1 0 0 0 0 0 0 0 0 0 0 0 0 0 
% 0 0 0 0 0 0 1 0 0 0 0 0 0 0 0 0 0 0 0 
% 0 0 0 0 0 0 0 1 0 0 0 0 0 0 0 0 0 0 0 
% 0 0 0 0 0 0 0 0 1 0 0 0 0 0 0 0 0 0 0 
\begin{verbatim}
 1 0 0 0 0 0 0 0 0 1 0 0 0 0 0 0 0 0 1 0 0 0 0 0 0 0 0
 0 1 0 0 0 0 0 0 0 0 1 0 0 0 0 0 0 0 0 1 0 0 0 0 0 0 0
 0 0 1 0 0 0 0 0 0 0 0 1 0 0 0 0 0 0 0 0 1 0 0 0 0 0 0
 0 0 0 1 0 0 0 0 0 0 0 0 1 0 0 0 0 0 0 0 0 1 0 0 0 0 0
 0 0 0 0 1 0 0 0 0 0 0 0 0 1 0 0 0 0 0 0 0 0 1 0 0 0 0
 0 0 0 0 0 1 0 0 0 0 0 0 0 0 1 0 0 0 0 0 0 0 0 1 0 0 0
 0 0 0 0 0 0 1 0 0 0 0 0 0 0 0 1 0 0 0 0 0 0 0 0 1 0 0
 0 0 0 0 0 0 0 1 0 0 0 0 0 0 0 0 1 0 0 0 0 0 0 0 0 1 0
 0 0 0 0 0 0 0 0 1 0 0 0 0 0 0 0 0 1 0 0 0 0 0 0 0 0 1
 1 0 0 1 0 0 1 0 0 0 0 0 0 0 0 0 0 0 0 0 0 0 0 0 0 0 0
 0 1 0 0 1 0 0 1 0 0 0 0 0 0 0 0 0 0 0 0 0 0 0 0 0 0 0
 0 0 1 0 0 1 0 0 1 0 0 0 0 0 0 0 0 0 0 0 0 0 0 0 0 0 0
 0 0 0 0 0 0 0 0 0 1 0 0 1 0 0 1 0 0 0 0 0 0 0 0 0 0 0
 0 0 0 0 0 0 0 0 0 0 1 0 0 1 0 0 1 0 0 0 0 0 0 0 0 0 0
 0 0 0 0 0 0 0 0 0 0 0 1 0 0 1 0 0 1 0 0 0 0 0 0 0 0 0
 1 1 1 0 0 0 0 0 0 0 0 0 0 0 0 0 0 0 0 0 0 0 0 0 0 0 0
 0 0 0 1 1 1 0 0 0 0 0 0 0 0 0 0 0 0 0 0 0 0 0 0 0 0 0
 0 0 0 0 0 0 0 0 0 1 1 1 0 0 0 0 0 0 0 0 0 0 0 0 0 0 0
 0 0 0 0 0 0 0 0 0 0 0 0 1 1 1 0 0 0 0 0 0 0 0 0 0 0 0
\end{verbatim}
Remember that the columns of the given array are the generators of the
semigroup.
The Hilbert basis of the cone generated by these $27$ vectors is the same as
the set of these $27$ vectors (we computed the Hilbert basis via 
{\tt normaliz}), so we conclude that the semigroup is saturated.
$3 \times 3 \times 5$ with $2$-marginals, 
after removing redundant rows (we removed redundant rows using {\tt cddlib}),
the constraint matrix becomes a $29 \times 45$ matrix.  Thus the semigroup
is generated by $45$ vectors in $\Z^{29}$. Note that these $45$ vectors are
extreme rays of the cone (we verified using {\tt cddlib}).  
The Hilbert basis of the cone of $A$ consists of these $45$ extreme rays 
so that we conclude that the semigroup is saturated.
Also $3 \times 3 \times 6$ with $2$-marginals, 
after removing redundant rows (we removed redundant rows using {\tt cddlib}),
the constraint matrix becomes a $34 \times 54$ matrix.  Thus the semigroup
is generated by $54$ vectors in $\Z^{34}$.
}

\section{Time complexity}\label{disc}

In 2002, \cite{Barvinok2002} introduced an algorithm to encode all 
integral vectors $\Bb \in \Z^d$ in Problem \ref{problem1}
as a {\em short rational generating function}
in polynomial time when $d$ and $n$ are fixed (Lemma \ref{genFunQ}
stated below). 
%In theory, Problem \ref{problem1} can be solved in polynomial time 
%in terms of the input size in fixed $d$ and $n$. 
However, in a sub-step of the
algorithm they use the {\em Projection Theorem} %(Lemma \ref{project})
which is not implementable
at present.  Thus, we do not know whether it is practical 
or not. From Lemma \ref{genFunQ}, % which is stated below, 
we can show that the time complexity of
computing $H$ is polynomial time if we fix $d$ and $n$ (Corollary \ref{cor:H}).

One might ask the time complexity of Problem \ref{problem2}. 
Using the results from [\cite{bar, BarviPom, Barvinok2002}], we can
prove  that Problem \ref{problem2} can be solved in 
polynomial time in fixed $d$ and $n$ (Theorem \ref{cputimeMain}). 
In order to prove the theorem, we will use the {\em multivariate generating 
function} of a set $X \subset \Z^d$, $f(X; x)$.
Namely, if $X\subset\Z^d$, define the generating function
\[f(X; x)=\sum_{s\in X}x^s,\]
where $x^s$ denotes $x_1^{s_1}\cdots x_d^{s_d}$ with
$s=(s_1,\ldots,s_d)$.   
If $X = P \cap \Z^d$ with fixed $d$, where $P$ is a rational convex polyhedron,  
or if $X = Q$ with fixed $d$ and $n$, then \cite{bar} and \cite{Barvinok2002}, respectively, showed that $f(X; x)$ can be written in the form of a 
polynomial-size sum of rational function of the form:
\begin{equation}\label{rational}
f(X; x) = \sum_{i \in I} \gamma_i \frac{x^{\alpha_i}}{\prod_{j = 1}^d (1 - x^{\beta_{ij}})}.
\end{equation}
Herein, $I$ is a finite (polynomial size) index set and all the
appearing data $\gamma_i\in\Q$ and $\alpha_i,\beta_{ij}\in\Z^d$ is
of size polynomial.
%As an example, if $P$ is the one-dimensional polytope $[0,N]$, $N \in \Z_+$,
%then $f(P \cap \Z;x)=1+x+x^2+\cdots+x^N$, $f(P \cap \Z;x)$ can be represented by a
%short rational function $(1-x^{N+1})/(1-x)$.
%We define such generating functions as {\em short rational functions} because
%it has a monomial in its numerator.
% Note that the rational function representation of $f(X; x)$ in
% \eqref{rational} is ``short.''  
% Therefore we say generating functions in 
% the form of \eqref{rational} as {\em short rational generating functions}.
If a rational generating function $f(X; x)$ is polynomial size in the total 
bit size of inputs, then $f(X; x)$ is called a {\em short rational generating 
function}.
As an example, if $P$ is the one-dimensional polytope $[0,N]$, $N \in \Z_+$,
then $f(P \cap \Z;x)=1+x+x^2+\cdots+x^N$, $f(P \cap \Z;x)$ can be represented 
by a short rational generating function $(1-x^{N+1})/(1-x)$.

\begin{thm}\label{cputimeMain}
Suppose we fix $d$ and $n$. %With given matrix $A$,
There is a polynomial time algorithm 
in terms of the input size
to decide whether the set of holes, $H$,
for the semigroup, $Q$, generated by the columns of $A$ is finite or not.
\end{thm}

Using the generating functions, we can show that the computation of fundamental
holes for $Q$ can be solved polynomial time if we fix $d$ and $n$. 
\begin{thm}\label{cputime}
Suppose we fix $d$ and $n$. Suppose $Q$ is not saturated.
The set of fundamental holes, $H_0$, can be  encoded in a short rational 
generating function 
in time polynomial in terms of the input size.
\end{thm}
One notes that this algorithm outputs a generating function in the form of a short rational generating function.  Therefore this does not return an explicit 
representation of $H_0$.  However, if one wants to enumerate all elements
in $H_0$, one can do the following:
from the proof of Proposition \ref{prop2}, we have $H_0 \subset P \cap
\Z^d$, where
\begin{equation}
\label{eq:compactP}
P := \{ x \in \R^d: x = \sum_{i=1}^n \delta_i \Ba_i, \, 0 \leq \delta_i
\leq 1\}.
\end{equation}
This shows that $H_0$ is finite and also gives a finite procedure to
enumerate $H_0$:
\begin{itemize}
\item Compute the Hilbert basis $B$ of $\cone(\Ba_1, \dots, \Ba_n)\cap L$.
\item Check each $z\in B$ whether it is a fundamental hole or not,
that is, compute $B\cap H_0$.
\item Generate all nonnegative integer combinations of elements in $B\cap
H_0$ that lie in $P \cap \Z^d$ and check for each such $z$ whether it is a
fundamental hole or not.
\end{itemize}
For more details, see [\cite{Hemmecke}].

Before proofs of Theorem \ref{cputimeMain} and Theorem \ref{cputime}, we would like to state lemmas from
[\cite{Barvinok2002}] and [\cite{BarviPom}].
\begin{lemma}[(7.3) in {[\cite{Barvinok2002}]}]\label{genFunQ}
Suppose we fix $d$ and $n$.  Let $Q = Q(A)$.  Then the generating function
$f(Q; x)$ for the semigroup $Q$ can be computed in polynomial time in terms
of the input size as a short rational generating function in the form of 
\eqref{rational}.
\end{lemma}
\begin{lemma}[Theorem 4.4 in {[\cite{BarviPom}]}]\label{genFunQsat}
Suppose we fix $d$ and suppose $P \subset \R^d$ is a rational convex 
polyhedron. Then the generating function $f(P \cap \Z^d; x)$ %for $P \cap \Z^d$
can be computed in polynomial time in terms
of the input size as a short rational generating function in the form of 
\eqref{rational}.
\end{lemma}
By Lemma \ref{genFunQ} and Lemma \ref{genFunQsat}, immediately, 
we have the following result.
\begin{cor}\label{cor:H}
Suppose we fix $d$ and $n$.  Let $Q = Q(A)$.  Then the generating function
$f(H; x)$ for the set of holes, $H := \Qsat \backslash Q$, can be computed in 
polynomial time in terms of the input size as a short rational generating
function in the form of \eqref{rational}.
\end{cor}
\begin{proof}
Suppose we fix $d$ and $n$. 
By Lemma \ref{genFunQ}, we can compute the generating function 
$f(Q; x)$ for the semigroup $Q$ in polynomial time and by Lemma 
\ref{genFunQsat} we can compute the generating function 
$f(\Qsat; x)$ for the semigroup $\Qsat$ in polynomial time.
The generating function $f(H; x)$ for $H$ is $f(\Qsat; x) - f(Q; x)$.
\end{proof}

Using Corollary \ref{cor:H}, %theorem \ref{cputimeMain}, 
we can prove Theorem \ref{cputimeMain}.
\begin{proof}[Proof of Theorem \ref{cputimeMain}]
Suppose we fix $d$ and $n$. 
First, we use Corollary \ref{cor:H} to compute the generating function, 
$f(H; x)$, for $H$ in polynomial time in the form of \eqref{rational}.
Let 
\[
f(H; x) = \sum_{i \in I} \gamma_i \frac{x^{\alpha_i}}{\prod_{j = 1}^d (1 - x^{\beta_{ij}})}.
\]
Then, we will do the following: First we choose $l \in \Z^d$ so that 
$\langle l, \beta_{ij}\rangle 
\neq 0$.  We find such $l$ in polynomial time by  Lemma 2.5 
in [\cite{Barvinok2002}].  Let $l=(\lambda_1, \dots, \lambda_d) \in 
\Z^d$.  For $\tau >  0$, let $x_{\tau} = (\exp(\tau\lambda_1), \dots, 
\exp(\tau\lambda_d))$ and let $\xi_{ij} = \langle l, \beta_{ij}\rangle $ and $\nu_i = 
\langle l, \alpha_i\rangle $.  Then we apply the monomial substitution
$x_i \to \exp(\tau \lambda_i)$.  We can do this monomial substitution 
in polynomial time by Lemma 2.5 and Theorem 2.6 in [\cite{Barvinok2002}].  %(5.2) in [\cite{BarviPom}].
Then
\[
f(H; x_{\tau}) = \frac{1}{\tau^d}\left(\sum_{i \in I} \gamma_i\frac{\tau^d\exp(\tau \nu_i)}{\prod_{j = 1}^d (1 - \exp(\tau \xi_{ij}))}\right).
\]
Let 
\[
h_i(\tau) = \frac{\tau^d\exp(\tau \nu_i)}{\prod_{j = 1}^d (1 - \exp(\tau \xi_{ij}))}
\]
is a holomorphic function in a neighborhood of $\tau = 0$ and we take the 
Taylor expansion around $\tau = 0$ (i.e., we take the Laurent expansion around
$\tau = 0$ for $h_i(\tau)/\tau^d$).
The coefficients of the $k$th powers, where $0 \leq k \leq d-1$, of the 
Taylor expansion of $h_i$ are:
\[
\frac{1}{\xi_{i1}\cdots \xi_{id}}\left(\sum_{l=0}^k \frac{\nu_i^l}{l!}
{\rm td}_{k-l}(\xi_{i1}, \dots,\xi_{id})\right),
\]
where ${\rm td}_l(\xi_{i1}, \dots,\xi_{id})$ is a homogeneous polynomial of
degree $l$ and which is called the $l$th {\em Todd polynomial} in $\xi_{i1}, 
\dots,\xi_{id}$ (see more details in Definition 5.1 in [\cite{BarviPom}]).

Now we claim that if the coefficients of negative powers of the Laurent 
expansion of 
$(\sum_{i \in I} h_i(\tau))/\tau^d$ 
are all canceled, then $H$ has to be finite.  We prove this by contradiction.
Suppose $H$ is infinite.  
Then, since all  coefficients of negative powers in the Laurent expansion
are canceled, the sum of the coefficients of the constant 
terms:
\begin{equation}\label{coef}
\sum_{i \in I} \frac{\gamma_i}{ \xi_{i1}\cdots \xi_{id}}
\left(\sum_{l=0}^d \frac{\nu_i^l}{l!}
{\rm td}_{d-l}(\xi_{i1}, \dots,\xi_{id})\right)
\end{equation}
must be equal to the number of elements in $H$ when we send $\tau \to 0$ 
((5.2) [\cite{BarviPom}]). 
Thus, the sum of the coefficients of the constant 
terms in \eqref{coef} must be equal to infinity. 
Since $I$ is a finite index set, a coefficient of
the constant term in the Laurent expansion of some rational function must 
be infinite.  However, the Todd polynomials are
polynomials in $\C$ so it is impossible.  Thus we reach a contradiction.

Conversely, it is obvious that
%Now we must show that 
if the coefficients of negative powers of the 
Laurent expansion of $(\sum_{i \in I} h_i(\tau))/\tau^d$ 
are not canceled, then $H$ is infinite.  
%This is obvious.

Therefore we will have to check all coefficients of the $k$th powers, 
where $0 \leq k \leq d-1$, of the Taylor expansion of each $h_i(\tau)$.
Since we have the polynomial size index set $I$ and we have to only check
$d$ coefficients for each rational function, this computation can be done in
polynomial time.
\end{proof}

Now we would like to discuss the {\em intersection algorithm}, which we
need to encode $H_0$ in a short rational 
generating function in polynomial time in fixed $d$ and $n$.
\begin{lemma}[Theorem 3.6 in \cite{Barvinok2002}] \label{intersect} 
Let $S_1, S_2$ be finite
subsets of $\Z^d$, for fixed $d$. Let $f(S_1;x)$ and $f(S_2;x)$ be their
 generating functions, given as short
rational generating functions with at most $k$ binomials in each denominator.
Then there exist a polynomial time algorithm, which, given
$f(S_i;x)$, computes 
$$ f(S_1 \cap S_2; x) \quad =  \quad \sum_{i \in I} \gamma_i \cdot \frac { x^{u_i} } {  (1-x^{v_{i1}})  \cdots (1-x^{v_{is}}) }$$
with $s \leq 2k$, where the $\gamma_i$ are rational numbers, $u_i,v_{ij}$ nonzero
integers, and $I$ is a polynomial-size index set.
\end{lemma}
The essential step in the {\em intersection algorithm} is
the {\em Hadamard product} [Definition 3.2 in \cite{Barvinok2002}].
Using Lemma \ref{intersect}, we can compute the union of $s$ sets in 
$\Z^d$ in polynomial time for fixed $d$ and $s$.

% \comment{
% Another key subroutine introduced by Barvinok and Woods is the
% following \emph {Projection Theorem}. In both Lemmas \ref{intersect}
% and \ref{project}, the dimension $d$ is assumed to be fixed.
% \begin{lemma}[Theorem 1.7 in \cite{Barvinok2002}] \label{project} 
% Assume the dimension $d$ is a fixed constant.
% Consider  a rational polytope $P \subset \R^d$ and a
% linear map $T: \Z^d \rightarrow \Z^k$.
% There is a polynomial time algorithm which
% computes a short representation of the
% generating function $\, f \bigl(T(P \cap \Z^d);x\bigr) $.
% \end{lemma}}

% Also from the proof of Proposition \ref{prop2}, we have $H_0 \subset P \cap
% \Z^d$, where
% \[
% P := \{ x \in \R^d: x = \sum_{i=1}^n \delta_i \Ba_i, \, 0 \leq \delta_i
% \leq 1\}.
% \]
% Also we have the following lemma.
% \begin{lemma}\label{fundHole}
% Let $P := \{ x \in \R^d: x = \sum_{i=1}^n \delta_i \Ba_i, \, 0 \leq \delta_i
% \leq 1\}$.  Then $H_0 \subset P \cap \Z^d$.
% \end{lemma}
% \begin{proof}
% A proof is in the proof of Proposition \ref{prop2}.
%By contradiction.  Suppose $y \in H_0$ can be written as $y = 
%\sum_{i=1}^n \delta_i \Ba_i$, with some $\delta_i > 1$.
%Note that $y - \Ba_i \in \Qsat$.  Let $y_0 = y - \Ba_i$.
%We claim that $y_0$ is a hole.  If not, then $y_0 \in Q$.  Then $y = y_0 +
%\Ba_i \in Q$.  Thus contradiction.  Therefore $y_0$ is a hole.
%But it contradicts that $y$ is fundamental.
%\end{proof}
We now give a proof of Theorem \ref{cputime}.

\begin{proof}[Proof of Theorem \ref{cputime}]
%Assume that we fix $n$ and $d$. 
Suppose $Q$ is not saturated.
% Let $P := \{ x \in \R^d: x = \sum_{i=1}^n \delta_i \Ba_i, \, 0 \leq \delta_i
% \leq 1\}$. %and let $P_0$ be the unit cube in $\R^n$. 
%Let $T: \Z^n \rightarrow \Z^d$ be 
%a linear transformation such that $x = \sum_{j=1}^n \xi_j \Ba_j$.
%Note that the semigroup $Q$ can be represented as the image of $T(\Z^n)$.
%Also note that $P \cap Q = T(P_0 \cap \Z^n)$.
Using Lemma \ref{genFunQsat}, we compute the generating function $f(P \cap \Z^d; x)$ in polynomial time, where $P$ is given in (\ref{eq:compactP}).
%Then, using Lemma \ref{project}, we compute the generating function
%$g(x)$ for $T(P_0\cap \Z^n) = P \cap Q$ in polynomial time.  
Note that there are $2^n$ points in $P \cap Q$, namely $\{x \in Q:
x = \sum_{i=1}^n \xi_i \Ba_i, \, \xi_i \in \{0, 1\}\}$.  So we can 
enumerate all points in $P \cap Q$ in constant time.
%By Lemma \ref{fundHole}, we know that 
%Fundamental holes are contained in
%$\bar H = P \cap \Z^d - P\cap Q = P \cap \Z^d - T(P_0 \cap \Z^n)$.  
Let $\bar H = (P \cap \Z^d) \setminus(P\cap Q)$. 
Its generating function $f(\bar H; x)$ is 
$f(P \cap \Z^d; x) - f(P \cap Q; x)$ and it can be computed in polynomial time.
Note that $H_0 = \bar H \setminus ((\bar H + (P \cap Q)) \cap \bar H)$ from 
the definition of $H_0$ and $H_0 \subset (P\cap \Z^d)\setminus(P \cap Q)$.

We compute the generating function for $(\bar H + (P \cap Q))$ by the 
following:  let $(P \cap Q) \setminus \{0\} = \{z_1, \cdots, z_{2^n-1}\}$.
For each $i = 1, 2, \dots, 2^n-1$, let $g_i(x) := x^{z_i} \cdot f(\bar H; x)$
which is the generating function for the set $z_i + \bar H$.
Since $2^n - 1$ is a constant (we are fixing $n$ as a constant), 
applying  Lemma \ref{intersect} we can compute the generating
function for the union of $z_i + \bar H$ in polynomial time. 
%The generating function for $(\bar H + (P \cap Q))$ is $g(x) \cdot f(\bar H; x)$.  
Since $\bar H$ and $(P \cap Q)$ are finite
%, we apply  Lemma \ref{intersect} and 
we are done.
\end{proof}

\noindent {\large\bf Acknowledgment}\quad 

We would like to thank Dr. Seth Sullivant for 
useful references,  % and conversation and to thank Raymond 
%Hemmecke and Matthias Walter for Hilbert bases computation via {\tt 4ti2}.
%Also we would like to thank 
Dr. Raymond Hemmecke and Prof. Robert Weismantel for
useful comments and references for Theorem \ref{condition}, Prof. Ezra
Miller for algebraic interpretations and proofs of our
results, Prof. Alexander Barvinok and Dr. Sven Verdoolaege 
for useful suggestions for the proof of 
Theorem \ref{cputimeMain},
Prof. Hidefumi Ohsugi for pointing out the relation between holes
and non-squarefree indispensable moves.  
At the last but not least,
we would like to thank Prof. Jesus De Loera for useful comments.

\par

\bibliographystyle{ims2}
\bibliography{satss}

\end{document}